\documentclass[11pt,reqno]{amsart}

\usepackage{amsmath,amsthm,amssymb}

\usepackage{tikz-cd} 

\usepackage{graphicx, psfrag, amsmath, amscd, amssymb,mathtools}
\usepackage{lipsum}

\usepackage[bookmarksnumbered, colorlinks, plainpages]{hyperref}
\hypersetup{colorlinks=true,linkcolor=red, anchorcolor=green, citecolor=cyan, urlcolor=red, filecolor=magenta, pdftoolbar=true}

\newtheorem{thm}{Theorem}[section]
\newtheorem{prop}[thm]{Proposition}
\newtheorem{lem}[thm]{Lemma}
\newtheorem{cor}[thm]{Corollary}

\theoremstyle{definition}
\newtheorem{defi}[thm]{Definition}
\newtheorem{exa}[thm]{Example}
\newtheorem{rem}[thm]{Remark}

\newcommand{\B}{{\mathcal{B}}}
\newcommand{\C}{{\mathcal{C}}}
\newcommand{\F}{{\mathcal{F}}}
\newcommand{\G}{{\mathscr{G}}}
\newcommand{\I}{{\mathscr{I}}}

\newcommand{\ch}{{\mathrm{ch}}}

\newcommand{\conv}{{\mathrm{conv}}}
\newcommand{\ind}{{\mathrm{ind}}}
\newcommand{\rk}{{\mathrm{rk}}}

\newcommand{\N}{{\mathcal{N}}}

\renewcommand{\I}{{\mathcal{I}}}
\renewcommand{\S}{{\mathfrak{S}}}

\newcommand{\x}{{\mathbf{x}}}
\renewcommand{\L}{{\mathcal{L}}}
\renewcommand{\G}{{\mathcal{G}}}

\begin{document}

\title[Stembridge codes and Chow rings]{Stembridge codes and Chow rings}

\author{Hsin-Chieh Liao}
\address{Department of Mathematics, University of Miami, USA}
\email{h.liao@math.miami.edu}

\date{November 15, 2022; revised December 10, 2022}

\maketitle

\begin{abstract}
It is well known that the Eulerian polynomial is the Hilbert series of the cohomology of the permutohedral variety. We answer a question of Stembridge on finding a geometric explanation of the \emph{permutation representation} this cohomology carries. Our explanation involves an $\mathfrak{S}_n$-equivariant bijection between a basis for the Chow ring of the Boolean matroid and codes introduced by Stembridge. There are analogous results for the stellohedral variety. We provide a geometric explanation of the permutation representation that its cohomology carries. This involves the augmented Chow ring of a matroid introduced by Braden, Huh, Matherne, Proudfoot and Wang. Along the way, we also obtain some new results on augmented Chow rings.    
\end{abstract}

\section{Introduction}
Consider the $(n-1)$-dimensional permutohedron 
\[
    \Pi_n\coloneqq\conv\{(\sigma_1,\ldots,\sigma_n)\in\mathbb{R}^n : \sigma=\sigma_1\sigma_2\ldots\sigma_n\in\S_n\}.
\] 
Its normal fan $\Sigma_n=\Sigma(\Pi_n)$ can be obtained from the \emph{braid arrangement} $H_{i,j}\coloneqq\{x\in\mathbb{R}^n : x_i=x_j\}$ for $1\le i<j\le n$. The toric variety $X_{\Sigma_n}$ associated to $\Sigma_n$ is called the \emph{permutohedral variety}. An intriguing fact about $X_{\Sigma_n}$ is that the Hilbert series of its cohomology $H^*(X_{\Sigma_n})$ is the Eulerian polynomial.

The cohomology $H^*(X_{\Sigma_n})$ carries a representation of $\S_n$ induced by the $\S_n$-action on $\Sigma_n$. Using work of Procesi  \cite{Procesi1990}, Stanley \cite{Stanley1989} computed its Frobenius series, which shows this representation is a permutation representation. 
Stembridge \cite{Stembridge1992} introduced a combinatorial object called a \emph{code} and showed that the representation of $\S_n$ on the space generated by codes has the same Frobenius series as $H^*(X_{\Sigma_n})$. He then asked if there is a geometric explanation of the permutation representation on $H^*(X_{\Sigma_n})$ (see \cite[p.317]{Stembridge1992}, \cite[p.252]{Stembridge1994}).

In this extended abstract\footnote{see \cite{liao2022} for the full paper}, we answer Stembridge's question by identifying $H^*(X_{\Sigma_n})$ with the Chow ring of the Boolean matroid and then finding a permutation basis for the induced action of $\S_n$ on the Chow ring. We show that a basis of Feichtner--Yuzvinsky for general matroids serves the purpose when we apply it to Boolean matroids. We do this by constructing an $\S_n$-equivariant bijection between this basis and Stembridge codes.

There is a parallel story involving the \emph{stellohedron} $\widetilde{\Pi}_n$. The story begins with the \emph{binomial Eulerian polynomial}, which Postnikov, Reiner, Williams \cite{PRW2008} show is equal to the Hilbert series of the cohomology of the toric variety associated to $\widetilde{\Pi}_n$. Shareshian and Wachs~\cite{ShareshianWachs2020} show that the representation of $\S_n$ on this cohomology is a permutation representation. We provide a geometric explanation for this, which involves the \emph{augmented Chow ring of matroids} introduced by Braden, Huh, Matherne, Proudfoot, Wang~\cite{BHMPW2020}, \cite{BHMPW2020+}. We introduce \emph{extended codes} in order to obtain results analogous to those for the permutohedron mentioned above. Along the way, we also obtain some general results for the augmented Chow rings and the augmented Bergman fans of matroids.


\section{Eulerian story: Permutohedra}

\subsection{The $\S_n$-module structure on $H^*(X_{\Sigma_n})$}


Given a $d$-dimensional simple polytope $P$ with normal fan $\Sigma(P)$ and dual polytope $P^*$, the $h$-polynomial $h_{P^*}(t)$ of $P^*$ agrees with the Hilbert series of the cohomology $H^*(X_{\Sigma(P)})$ of the toric variety $X_{\Sigma(P)}$ (see \cite{Stanley1989} eq. (26)). 
It is well-known that the $h$-polynomial $h_{\Pi_n^*}(t)$ for the dual permutohedron $\Pi_n^*$ is the Eulerian polynomial $A_n(t)$, hence one has
\[
    A_n(t)=h_{\Pi_n^*}(t)=\sum_{j=0}^{n-1}\dim H^{2j}(X_{\Sigma_n})t^j.
\]
The cohomology $H^*(X_{\Sigma_n})$ carries an $\mathfrak{S}_n$-representation induced by the $\S_n$-action on $\Sigma_n$. Using the work of Procesi \cite{Procesi1990}, Stanley \cite{Stanley1989} computed the Frobenius series of this representation:
\begin{equation} \label{FrobPerm}
    \sum_{n\ge 0}\sum_{j=0}^{n-1}\ch\left(H^{2j}(X_{\Sigma_n})\right) t^jz^n=\frac{(1-t)H(z)}{H(tz)-tH(z)},
\end{equation}
where $\ch$ is the \emph{Frobenius characteristic}, $H(z)=\sum_{n\ge 0}h_n(\mathbf{x})z^n$ and $h_n$ is the complete homogeneous symmetric function of degree $n$. From (\ref{FrobPerm}), one can see that $H^*(X_{\Sigma_n})$ carries a permutation representation of $\S_n$. 

A \emph{Stembridge code} is a sequence $\alpha$ over $\{0,1,2,\ldots\}$ with marks such that if $m(\alpha)$ is the maximum number appearing in $\alpha$ then for each $k\in [m(\alpha)]$
\begin{itemize}
\item $k$ occurs at least twice in $\alpha$;
\item a mark is assigned to the $i$th occurence of $k$ for a unique $i\ge 2$.
\end{itemize}
For all $k\in [m(\alpha)]$, let $f(k)$ be the number of occurrences of $k$ in $\alpha$ to the left of the marked $k$. (So $f(k)=i-1$.) We let $(\alpha,f)$ denote the Stembridge code. The \emph{index} of  $(\alpha,f)$ is \[\ind(\alpha,f)\coloneqq\sum_{k\in [m(\alpha)]}f(k).\] 
Let $\C_n=\cup_{j=0}^{n-1}\C_{n,j}$ where $\C_{n,j}$ is the set of Stembridge codes of length $n$ with index $j$.
\begin{exa}
A Stembridge code $(\alpha,f)=11320\hat{2}\hat{3}\hat{1}2$ consists of $\alpha=11320231$ with $f(1)=2$, $f(2)=1$, $f(3)=1$ and $\ind(\alpha,f)=4$.
There are $6$ codes in $\C_3$: 
\[
\begin{array}{c|cccccc}
  (\alpha,f)    & 000 & 01\hat{1}&~10\hat{1}&~1\hat{1}0& 1\hat{1}1 & 11\hat{1}\\
  \hline
 \ind(\alpha,f) &   0  &    1     &     1    &    1     &    1    &  2                
\end{array}
\]
\end{exa}
For $\sigma\in\mathfrak{S}_n$, define $\sigma\cdot (\alpha_1\alpha_2,\ldots\alpha_n, f)=(\alpha_{\sigma(1)}\alpha_{\sigma(2)}\ldots\alpha_{\sigma(n)},f)$.
The action induces a graded representation $V_n=\bigoplus_{j=0}^{n-1}V_{n,j}$ of $\mathfrak{S}_n$, where $V_{n,j}=\mathbb{C}\mathcal{C}_{n,j}$. Its graded Frobenius series was computed by Stembridge and it agrees with (\ref{FrobPerm}). That is
\begin{equation}\label{CodeCohoQ}
	Q_n(\x,t)\coloneqq\sum_{j=0}^{n-1}\ch(V_{n,j})t^j=\sum_{j=0}^{n-1}\ch\left(H^{2j}(X_{\Sigma_n})\right)t^j.
\end{equation} 
Therefore $V_{n}\cong_{\S_n}H^*(X_{\Sigma_n})$. Stembridge \cite[p.317]{Stembridge1992} then asked for a geometric explanation of $H^*(X_{\Sigma_n})$ being a permutation representation, or more explicitly in \cite[p.252]{Stembridge1994}, whether there is a basis of $H^*(X_{\Sigma_n})$ permuted by $\S_n$ that induces the representation we are looking at. Although there is no obvious direct connection between Stembridge codes and $H^*(X_{\Sigma_n})$, it is natural that we expect such basis to have similar combinatorial structure as Stembridge codes. We will present such a basis in section \ref{CodeBasis}.

\subsection{Building sets and Chow rings of atomic lattices} \label{BS}
Here we recall the background about building sets, nested set complexes, and Chow rings of the atomic lattice from \cite{FeichtnerYuzvinsky2004}. For any poset $P$ and $X\in P$, we write $P_{\le X}=\{Y\in P: Y\le X\}$.


Let $\L$ be an atomic lattice. A subset $\mathcal{G}\subseteq \mathcal{L}-\{\widehat{0}\}$ is a \emph{building set} of $\mathcal{L}$ if for any $X\in\L-\{\widehat{0}\}$ with subset $\max(\mathcal{G}_{\le X})=\{G_1,\ldots,G_k\}$, there is a poset isomorphism 
\[
   \varphi_X: \prod_{i=1}^k[\hat{0},G_i]\longrightarrow [\hat{0},X]
\]
with $\varphi_X(\widehat{0},\ldots,\widehat{0},G_i,\widehat{0},\ldots,\widehat{0})=G_i$ for $i=1,\ldots,k$.

Given a building set $\G$ of $\L$, we say a subset $N\subseteq\G$ is \emph{nested} (or \emph{$\G$-nested}) if for any pairwise incomparable elements $G_1,\ldots, G_t\in N$ ($t\ge 2$), their join $G_1\vee\ldots\vee G_t\notin~\G$. Notice that the collection of all $\G$-nested sets forms an abstract simplicial complex $\mathcal{N}(\L,\G)$ which is called the \emph{nested set complex}. If $\G$ contains the maximal element $\widehat{1}$ then $\mathcal{N}(\L,\G)$ is a cone with apex $\{\widehat{1}\}$, in this case the base of the cone is called the \emph{reduced nested set complex} $\widetilde{\mathcal{N}}(\L,\G)$.

The nested set complex $\mathcal{N}(\L,\G)$ can be viewed as a generalization of the order complex $\Delta(\L-\{\widehat{0}\})$. Indeed, $\L-\{\widehat{0}\}$ is the maximal building set, $\mathcal{N}(\L,\L-\{\widehat{0}\})$ is $\Delta(\L-\{\widehat{0}\})$ and  $\widetilde{\mathcal{N}}(\L,\L-\{\widehat{0}\})$ is $ \Delta(\L-\{\widehat{0},\widehat{1}\})$.

\begin{defi}[\cite{FeichtnerYuzvinsky2004}]\label{ChowRingLattice}
Let $\L$ be an atomic lattice, $\mathfrak{A}(L)$ be the set of atoms and $\G$ be a building set of $\L$. Then the \emph{Chow ring of $\L$ with respect to $\G$} is the $\mathbb{Q}$-algebra
\[
	D(\L,\G)\coloneqq \mathbb{Q}[x_{G}]_{G\in\G}/(I+J)
\] 
where 
	$I=\left\langle\prod_{i=1}^t G_i: \{G_1,\ldots,G_t\}\notin\mathcal{N}(\L,\G)\right\rangle$  and $J=\left\langle\sum_{G\ge H}x_G: H\in\mathfrak{A}(\L)\right\rangle$.

\end{defi}
In particular, the \emph{Chow ring of a matroid} $M$ as defined in \cite{AHK2018}, \cite{BHMPW2020}, \cite{BHMPW2020+} is the special case that $\L$ is the lattice of flats of $M$ and $\G$ is the maximal building set, i.e. $D(\L(M),\L(M)-\{\emptyset\})$. Note that $\mathbb{Q}[x_{G}]_{G\in\G}/I$ is the Stanley-Reisner ring of $\mathcal{N}(\L,\G)$. For basic notions of matroid theory and Stanley-Reisner rings we refer the readers to \cite{Oxley2006matroid} and \cite{SturmfelsCommtative2005} respectively.

Feichtner and Yuzvinsky \cite{FeichtnerYuzvinsky2004} also found a Gr\"{o}bner basis for the ideal $I+J$ which gives the following basis for $D(\L,\G)$.

\begin{prop}[\cite{FeichtnerYuzvinsky2004}]\label{FYbasis} The following monomials form a basis for $D(\L,\G)$
\[
	\left\{\prod_{G\in N}x_G^{a_G}: N\text{ is nested, }a_G<d(G',G)\right\}
\]
where $G'$ is the join of $N\cap\L_{<G}$ and $d(G',G)$ is the minimal number $d$ such that $H_1,\ldots, H_d\in\frak{A}(\L)$ satisfies $G'\vee H_1\vee\ldots\vee H_d=G$. In particular, if $\L$ is a geometric lattice, then $d(G',G)=\rk(G)-\rk(G')$.
\end{prop}

\subsection{Stembridge codes and FY-basis for Chow ring of $B_n$} \label{CodeBasis}
In this section, we shall answer Stembridge's question. We identify $H^*(X_{\Sigma_n})$ with the Chow ring of the Boolean matroid, and use the facts in section \ref{BS} to find a basis permuted by the $\mathfrak{S}_n$-action.

The following lemma from \cite[p.251(1.5)]{Stembridge1994} allows us to identify  $H^*(X_{\Sigma_n})$ with the quotient of Stanley-Reisner ring of the boundary $\partial\Pi_n^*$ of the dual permutohedron. 
\begin{lem}[\cite{danilov1978},\cite{Stembridge1994}]\label{CohomologyFaceRing}
Let $P$ be a simple $n$-lattice polytope in an Euclidean $n$-space $V$ with normal fan $\Sigma(P)$ or equivalently the face fan of $P^*$. Let $K[\partial P^*]$ be the Stanley-Reisner ring of $\partial P^*$ over a field $K$ of characteristic $0$. Let $\theta_i=\sum_{v\in V(P^*)}\langle v,e_i \rangle x_v$ for $i=1,\ldots,n$, where $\{e_i\}$ is the standard basis in $\mathbb{Z}^n\subset V$, and $V(P^*)$ is the set of vertices of $P^*$. If a finite group $G$ acts on $\Sigma(P)$ simplicially and freely, then 
\[
	H^*(X_{\Sigma(P)},K)\cong K[\partial P^*]/\langle\theta_1,\ldots,\theta_n\rangle ~~\text{ as $K[G]$-modules}.
\]
\end{lem}

Let $M$ be a loopless matroid on ground set $[n]$ with a lattice of flats $\L(M)$. For each $S\subseteq [n]$, write $e_S=\sum_{i\in S}e_i$. The \emph{Bergman fan} $\Sigma_M$ of $M$ is a fan in  $\mathbb{R}^n/\langle e_{[n]}\rangle$ consists of cones $\sigma_{\F}$ indexed by all flags $\F=\{F_1\subsetneq\ldots\subsetneq F_k\}$ in $\L(M)-\{\emptyset,[n]\}$, where
\[
	\sigma_{\F}=\mathbb{R}_{\ge 0}\{e_{F_1},\ldots,e_{F_k}\}.
\]
The \emph{Bergman complex} is the simplicial complex obtained by intersecting $\Sigma_M$ with the unit sphere centered at $0$. Obviously, the Bergman complex is a geometric realization of the order complex $\Delta(\L(M)-\{\emptyset,[n]\})$.
 
The \emph{Chow ring of $M$} has the following two presentations
\begin{align}
    A(M) &\coloneqq\frac{\mathbb{Q}[x_F]_{F\in\mathcal{\L}(M)-\{\emptyset\}}/\left\langle x_F x_G: F,G\text{ are incomparable in }\mathcal{L}(M)\right\rangle}{\left\langle\sum_{F: i\in F}x_F: 1\le i \le n \right\rangle} \label{FYPresentation}\\
    \nonumber\\
    	&=\frac{\mathbb{Q}[x_F]_{F\in\mathcal{\L}(M)-\{\emptyset,[n]\}}/\left\langle x_F x_G: F,G\text{ are incomparable in }\mathcal{L}(M)\right\rangle}{\left\langle\sum_{F: i\in F}x_F-\sum_{F: j\in F}x_F: i\neq j \right\rangle}. \label{HuhPresentation}
\end{align}
The presentation (\ref{FYPresentation}) is a special case of Definition \ref{ChowRingLattice}, and (\ref{HuhPresentation}) is obtained from (\ref{FYPresentation}) by eliminating $x_E$. The presentation (\ref{HuhPresentation}) are used in \cite{AHK2018},\cite{BHMPW2020},\cite{BHMPW2020+}. Note that the numerator in (\ref{HuhPresentation}) is the Stanley-Reisner ring of $\Sigma_M$, or equivalently the Bergman complex of $M$.

Now let $M$ be the Boolean matroid $B_n$; the flats of $B_n$ are all subsets of $[n]$ and $\L(B_n)$ is the Boolean lattice. It is well-known that the order complex $\Delta(\L(B_n)-\{\emptyset,[n]\})$ is the  $\partial\Pi_n^*$. After some easy calculations, applying Lemma \ref{CohomologyFaceRing} to (\ref{HuhPresentation}) shows that
\[
    A(B_n)=\frac{\mathbb{Q}[\partial\Pi_n^*]}{\langle\theta_1,\ldots,\theta_{n-1}\rangle}\cong_{\mathfrak{S}_n} H^*(X_{\Sigma_n},\mathbb{Q}).
\]
By applying Proposition \ref{FYbasis} to (\ref{FYPresentation}), the Feichtner--Yuzvinsky basis of $A(B_n)$ is given by
\[
    FY(B_n)\coloneqq\left\{x_{F_1}^{a_1}x_{F_2}^{a_2}\ldots x_{F_k}^{a_k}:\substack{~\emptyset=F_0\subsetneq F_1\subsetneq F_2\subsetneq\ldots\subsetneq F_k\subseteq [n],\\
    1\le a_i\le |F_i|-|F_{i-1}|-1}\right\}.
\]
Note that $|F_i|-|F_{i-1}|\ge 2$ for all $i$.
We see that $\mathfrak{S}_n$ permutes $FY(B_n)$ and makes $A(B_n)$ an $\mathfrak{S}_n$-module. It turns out $FY(B_n)$ has similar structure as the codes $\mathcal{C}_n$. 
\begin{thm}\label{BasisBijection1}
There is a bijection $\phi:FY(B_n)\rightarrow\mathcal{C}_n$ that respects the $\mathfrak{S}_n$-actions and takes the degree of the monomials to the index of the corresponding codes. 
\end{thm}

The bijection is defined as follows. Given $u=x_{F_1}^{a_1}\ldots x_{F_k}^{a_k}\in FY(B_n)$, let $\phi(u)=(\alpha,f)$, where 
$\alpha_i=\begin{cases} j & \text{ if }i\in F_j-F_{j-1}\\
0 & \text{ if }i\in [n]-F_k
\end{cases}$ for $i\in [n]$ ~and~ $f(j)=a_j$ for $j\in[k]$.\\ 

\begin{exa}
The basis for $A^2(B_4)$ and the corresponding codes in $\C_{4,2}$:\\
\begin{tikzcd}[ampersand replacement=\&, column sep=0]
x_{12}x_{1234} \& x_{13}x_{1234} \& x_{14}x_{1234} \& x_{23}x_{1234} \& x_{24}x_{1234} \& x_{34}x_{1234}
\& x_{123}^2 \& x_{124}^2 \& x_{134}^2 \&  x_{234}^2 \& x_{1234}^2
\\
1\hat{1}2\hat{2} \arrow[u,leftrightarrow] 
\& 12\hat{1}\hat{2} \arrow[u,leftrightarrow] 
\& 12\hat{2}\hat{1} \arrow[u,leftrightarrow] 
\& 21\hat{1}\hat{2} \arrow[u,leftrightarrow] 
\& 21\hat{2}\hat{1} \arrow[u,leftrightarrow] 
\& 2\hat{2}1\hat{1} \arrow[u,leftrightarrow]
\& 11\hat{1}0 		\arrow[u,leftrightarrow] 
\& 110\hat{1} \arrow[u,leftrightarrow] 
\& 101\hat{1} \arrow[u,leftrightarrow] 
\& 011\hat{1} \arrow[u,leftrightarrow] 
\& 11\hat{1}1 \arrow[u,leftrightarrow]
\end{tikzcd}
It is easy to see that the bijection respects the $\mathfrak{S}_4$-action on both sets. 
\end{exa}
Theorem \ref{BasisBijection1} shows that $A^j(B_n)\otimes\mathbb{C}$ and $V_{n,j}$ are isomorphic permutation representation for all $0\le j\le n-1$. To our knowledge, this is the first found basis that bears an explicit resemblance to codes. 

\section{Binomial Eulerian story: stellohedra}
We will now switch track to a parallel story. The \emph{binomial Eulerian polynomial} $\widetilde{A}_n(t)\coloneqq 1+t\sum_{k=1}^n\binom{n}{k}A_k(t)$ shares many similar combinatorial and geometric properties with the Eulerian polynomial (see Shareshian and Wachs \cite{ShareshianWachs2020}). Geometrically, it is shown by Postinkov, Reiner, Williams \cite{PRW2008} that $\widetilde{A}_n(t)$ is the $h$-polynomial of the dual of the \emph{stellohedron} $\widetilde{\Pi}_n$. Let $\Delta_n=\conv(0,e_i:i\in[n])$ where $e_i$ is the standard basis vector in $\mathbb{R}^n$. The \emph{stellohedron} $\widetilde{\Pi}_n$ is the simple polytope obtained from $\Delta_n$ by truncating all faces not containing $0$ in the inclusion order. Let $\widetilde{\Sigma}_n$ be its normal fan. The Hilbert series of  $H^*(X_{\widetilde{\Sigma}_n})$ satisfies
\[
	\widetilde{A}_n(t)=h_{\widetilde{\Pi}_n^*}(t)=\sum_{j\ge 0}\dim H^{2j}(X_{\widetilde{\Sigma}_n})t^j.
\]

Shareshian and Wachs \cite{ShareshianWachs2020} introduced the symmetric function analogue of $\widetilde{A}_n(t)$,
\begin{equation} \label{defBinomQ}
	\widetilde{Q}_n(\x,t)\coloneqq h_n(\x)+t\sum_{k=1}^n h_{n-k}(\x)Q_k(\x,t)
\end{equation}
and showed that if we consider the simiplicial action of $\mathfrak{S}_n$ on $\widetilde{\Pi}_n^*$, the induced representation of $\S_n$ on  $H^*(X_{\widetilde{\Sigma}_n})$ has the following graded Frobenius series.
\begin{thm}[\cite{ShareshianWachs2020}]\label{FrobStell}
For all $n\ge 1$, we have $\sum_{j=0}^{n}\ch\left(H^{2j}(X_{\widetilde{\Sigma}_n})\right)t^j=\widetilde{Q}_n(\x,t)$.
\end{thm} 
From (\ref{defBinomQ}), since $Q_k(\x,t)$ is $h$-positive, so is $\widetilde{Q}_n(\x,t)$. Thus $H^*(X_{\widetilde{\Sigma}_n})$ also carries a permutation representation of $\mathfrak{S}_n$.

\subsection{Extended codes} \label{Sec:ExtCodes}
We introduce an analog of Stembridge codes called the \emph{extended codes}.
 
 An \emph{extended code}  $(\alpha,f)$ is a marked sequence like a Stembridge code. The sequence $\alpha$ is over $\{0,1,\ldots\}\cup\{\infty\}$ with $\infty$'s working as $0$'s in codes; $m(\alpha)$ and $f$ are defined as in Stembridge codes. Define the \emph{index} of the extended code $(\alpha, f)$ as 
\[
    \ind(\alpha,f)=\begin{cases}
        -1 & \text{ if }\alpha=\infty\ldots\infty;\\
        \sum_{k\in m(\alpha)}f(k) & \text{ otherwise.}
    \end{cases}
\]
Let $\widetilde{\mathcal{C}}_{n,j}$ be the set of extended codes of length $n$ with index $j$ and $\widetilde{\mathcal{C}}_n=\bigcup_{j=-1}^{n-1}\widetilde{\mathcal{C}}_{n,j}$.
\begin{exa}
The following are all the extended codes of length $3$:\\
$\widetilde{\mathcal{C}}_{3,-1}=\{\infty\infty\infty\}$,~ 
$\widetilde{\mathcal{C}}_{3,0}=\{0\infty\infty  ,\infty 0\infty ,\infty \infty 0,\infty 00,0\infty0 ,00\infty ,000\}$,~\\
$\widetilde{\mathcal{C}}_{3,1}=\{1\hat{1}\infty,1\infty\hat{1},\infty 1\hat{1},01\hat{1},10\hat{1},1\hat{1}0, 1\hat{1}1\}$,~
$\widetilde{\mathcal{C}}_{3,2}=\{11\hat{1}\}$.
\end{exa}
For $\sigma\in\mathfrak{S}_n$, define $\sigma\cdot (\alpha_1\alpha_2,\ldots\alpha_n, f)=(\alpha_{\sigma(1)}\alpha_{\sigma(2)}\ldots\alpha_{\sigma(n)},f)$. This induces an graded $\mathfrak{S}_n$-representation on $\widetilde{V}_n=\bigoplus_{j=0}^{n}\widetilde{V}_{n,j-1}$ where $\widetilde{V}_{n,j-1}=\mathbb{C}\widetilde{\C}_{n,j-1}$. We compute its Frobenius series and obtain a result parallel to Stembrige codes.
\begin{thm}.
For $n\ge 1$, we have  $\sum_{j=0}^n\ch(\widetilde{V}_{n,j})t^j =\widetilde{Q}_n(\x,t)$.
\end{thm}
Combining with Theorem \ref{FrobStell}, we see $\widetilde{V}_n\cong_{\S_n} H^*(X_{\widetilde{\Sigma}_n})$ as permutation modules. One can also ask if there is a permutation basis for $H^*(X_{\widetilde{\Sigma}_n})$ which has similar combinatorial structure as extended codes. We will show such a basis exists in what follows.

\subsection{Augmented Chow ring of $B_n$ and $H^*(X_{\widetilde{\Sigma}_n})$} \label{Sec:AugChowBn}
Braden, Huh, Matherne, Proudfoot, and Wang \cite{BHMPW2020} recently introduced the \emph{augmented Bergman fan} and the \emph{augmented Chow ring} of a matroid. They showed that the augmented Bergman fan of $B_n$ is the normal fan $\widetilde{\Sigma}_n$ of the stellohedron $\widetilde{\Pi}_n$. Therefore the corresponding spherical complex is the boundary complex $\partial\widetilde{\Pi}_n^*$. Below we use Lemma~ \ref{CohomologyFaceRing} to identify $H^*(X_{\widetilde{\Sigma}_n})$ with the augmented Chow ring of $B_n$.

Let $M$ be a loopless matroid on $[n]$ with lattice of flats $\L(M)$ and the collection of independent subsets $\I(M)$. The \emph{augmented Chow ring of $M$} is defined as 
\begin{equation}\label{AugChow}
 \widetilde{A}(M):=\frac{\mathbb{Q}\left[\{x_F\}_{F\in\L(M)\setminus [n]}\cup\{y_1,y_2,\ldots,y_n\}\right]/(I_1+I_2)}{\langle y_i-\sum_{F:i\notin F}x_F\rangle_{i=1,2,\ldots,n}}
\end{equation}
where $I_1=\left\langle x_F x_G: F,G\text{ are incomparable in }\mathcal{L}(M)\right\rangle$, $I_2=\left\langle y_ix_F: i\notin F\right\rangle$. There is a simplicial fan associated with $\widetilde{A}(M)$ called the \emph{augmented Bergman fan of M}.\\
Let $I\in\I(M)$ and $\F=(F_1\subsetneq\ldots\subsetneq F_k)$ be a flag in $\L(M)$. We say $I$ is \emph{compatible with} $\F$, denoted by $I\le \F$, if $I\subseteq F_1$. Note that for the empty flag $\emptyset$, we have $I\le \emptyset$ for any $I\in\I(M)$. The \emph{augmented Bergman fan $\widetilde{\Sigma}_M$ of $M$} is a simplicial fan in $\mathbb{R}^n$ consisting of cones $\sigma_{I\le\F}$ indexed by compatible pairs $I\le\F$, where $\F$ is a flag in $\L(M)-\{[n]\}$ and
\[
	\sigma_{I\le\F}=\mathbb{R}_{\ge 0}\left(\{e_i\}_{i\in I}\cup\{-e_{[n]\setminus F}\}_{F\in\mathcal{F}}\right).
\]
The corresponding simplicial complex is called the \emph{augmented Bergman complex}.

\begin{exa} \label{ExAugBergB2}
The augmented Bergman fan $\widetilde{\Sigma}_{B_2}$ of the Boolean matroid $B_2$ in $\mathbb{R}^2$ and the corresponding augmented Bergman complex:

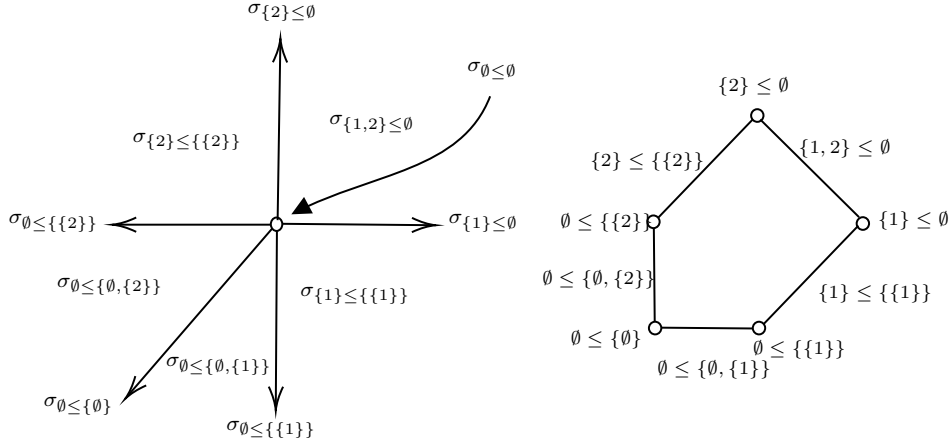
\tikzset{every picture/.style={line width=0.75pt}} 

\begin{center}
\begin{tikzpicture}[x=0.75pt,y=0.75pt,yscale=-1.1,xscale=1.1]

\draw    (130.09,119.65) -- (198.46,120.29) ;
\draw [shift={(200.46,120.31)}, rotate = 180.54] [color={rgb, 255:red, 0; green, 0; blue, 0 }  ][line width=0.75]    (10.93,-3.29) .. controls (6.95,-1.4) and (3.31,-0.3) .. (0,0) .. controls (3.31,0.3) and (6.95,1.4) .. (10.93,3.29)   ;
\draw   (124.84,120.05) .. controls (124.74,118.35) and (125.84,116.88) .. (127.29,116.77) .. controls (128.74,116.66) and (129.99,117.95) .. (130.09,119.65) .. controls (130.18,121.34) and (129.08,122.81) .. (127.63,122.92) .. controls (126.18,123.03) and (124.93,121.74) .. (124.84,120.05) -- cycle ;
\draw    (127.95,117.98) -- (129.2,36.63) ;
\draw [shift={(129.23,34.64)}, rotate = 90.88] [color={rgb, 255:red, 0; green, 0; blue, 0 }  ][line width=0.75]    (10.93,-3.29) .. controls (6.95,-1.4) and (3.31,-0.3) .. (0,0) .. controls (3.31,0.3) and (6.95,1.4) .. (10.93,3.29)   ;
\draw    (124.84,120.05) -- (54.18,120.05) ;
\draw [shift={(52.18,120.05)}, rotate = 360] [color={rgb, 255:red, 0; green, 0; blue, 0 }  ][line width=0.75]    (10.93,-3.29) .. controls (6.95,-1.4) and (3.31,-0.3) .. (0,0) .. controls (3.31,0.3) and (6.95,1.4) .. (10.93,3.29)   ;
\draw    (127.63,122.92) -- (127.11,203.32) ;
\draw [shift={(127.1,205.32)}, rotate = 270.37] [color={rgb, 255:red, 0; green, 0; blue, 0 }  ][line width=0.75]    (10.93,-3.29) .. controls (6.95,-1.4) and (3.31,-0.3) .. (0,0) .. controls (3.31,0.3) and (6.95,1.4) .. (10.93,3.29)   ;
\draw    (125.55,121.71) -- (59.36,198.44) ;
\draw [shift={(58.06,199.95)}, rotate = 310.78] [color={rgb, 255:red, 0; green, 0; blue, 0 }  ][line width=0.75]    (10.93,-3.29) .. controls (6.95,-1.4) and (3.31,-0.3) .. (0,0) .. controls (3.31,0.3) and (6.95,1.4) .. (10.93,3.29)   ;
\draw    (225.39,61.14) .. controls (212.13,95.44) and (171.32,96.33) .. (136.35,114.2) ;
\draw [shift={(134.22,115.31)}, rotate = 331.84] [fill={rgb, 255:red, 0; green, 0; blue, 0 }  ][line width=0.08]  [draw opacity=0] (8.93,-4.29) -- (0,0) -- (8.93,4.29) -- cycle    ;
\draw   (345.05,70.15) .. controls (344.95,68.55) and (346.13,67.17) .. (347.7,67.06) .. controls (349.26,66.96) and (350.61,68.17) .. (350.71,69.78) .. controls (350.81,71.38) and (349.63,72.76) .. (348.06,72.86) .. controls (346.5,72.97) and (345.15,71.75) .. (345.05,70.15) -- cycle ;
\draw   (393.44,119.7) .. controls (393.34,118.1) and (394.53,116.72) .. (396.09,116.62) .. controls (397.65,116.51) and (399,117.73) .. (399.11,119.33) .. controls (399.21,120.93) and (398.02,122.31) .. (396.46,122.42) .. controls (394.89,122.52) and (393.54,121.31) .. (393.44,119.7) -- cycle ;
\draw   (297.42,118.92) .. controls (297.32,117.32) and (298.51,115.93) .. (300.07,115.83) .. controls (301.64,115.73) and (302.99,116.94) .. (303.09,118.54) .. controls (303.19,120.14) and (302,121.52) .. (300.44,121.63) .. controls (298.88,121.73) and (297.53,120.52) .. (297.42,118.92) -- cycle ;
\draw   (298.19,167.68) .. controls (298.09,166.08) and (299.28,164.7) .. (300.84,164.59) .. controls (302.4,164.49) and (303.75,165.7) .. (303.86,167.31) .. controls (303.96,168.91) and (302.77,170.29) .. (301.21,170.39) .. controls (299.64,170.5) and (298.29,169.28) .. (298.19,167.68) -- cycle ;
\draw   (345.82,167.68) .. controls (345.72,166.08) and (346.9,164.7) .. (348.47,164.59) .. controls (350.03,164.49) and (351.38,165.7) .. (351.48,167.31) .. controls (351.58,168.91) and (350.4,170.29) .. (348.83,170.39) .. controls (347.27,170.5) and (345.92,169.28) .. (345.82,167.68) -- cycle ;
\draw    (345.82,71.73) -- (302.32,116.97) ;
\draw    (394.21,120.49) -- (350.71,165.73) ;
\draw    (349.23,71.78) -- (394.55,116.62) ;
\draw    (300.44,121.63) -- (300.84,164.59) ;
\draw    (345.82,167.68) -- (303.86,167.31) ;

\draw (205.12,113.9) node [anchor=north west][inner sep=0.75pt]  [font=\scriptsize] [align=left] {$\displaystyle \sigma _{\{1\} \leq \emptyset }$};
\draw (112.94,17.38) node [anchor=north west][inner sep=0.75pt]  [font=\scriptsize] [align=left] {$\displaystyle \sigma _{\{2\} \leq \emptyset }$};
\draw (103.09,208.74) node [anchor=north west][inner sep=0.75pt]  [font=\scriptsize] [align=left] {$\displaystyle \sigma _{\emptyset \leq \{\{1\}\}}$};
\draw (3.37,113.9) node [anchor=north west][inner sep=0.75pt]  [font=\scriptsize] [align=left] {$\displaystyle \sigma _{\emptyset \leq \{\{2\}\}}$};
\draw (19.33,197.74) node [anchor=north west][inner sep=0.75pt]  [font=\scriptsize] [align=left] {$\displaystyle \sigma _{\emptyset \leq \{\emptyset \}}$};
\draw (150.4,68.39) node [anchor=north west][inner sep=0.75pt]  [font=\scriptsize] [align=left] {$\displaystyle \sigma _{\{1,2\} \leq \emptyset }$};
\draw (213.66,44.05) node [anchor=north west][inner sep=0.75pt]  [font=\scriptsize] [align=left] {$\displaystyle \sigma _{\emptyset \leq \emptyset }$};
\draw (137.13,147.23) node [anchor=north west][inner sep=0.75pt]  [font=\scriptsize] [align=left] {$\displaystyle \sigma _{\{1\} \leq \{\{1\}\}}$};
\draw (60.2,76.56) node [anchor=north west][inner sep=0.75pt]  [font=\scriptsize] [align=left] {$\displaystyle \sigma _{\{2\} \leq \{\{2\}\}}$};
\draw (25.57,143.07) node [anchor=north west][inner sep=0.75pt]  [font=\scriptsize] [align=left] {$\displaystyle \sigma _{\emptyset \leq \{\emptyset ,\{2\}\}}$};
\draw (75.41,179.57) node [anchor=north west][inner sep=0.75pt]  [font=\scriptsize] [align=left] {$\displaystyle \sigma _{\emptyset \leq \{\emptyset ,\{1\}\}}$};
\draw (328.43,51.19) node [anchor=north west][inner sep=0.75pt]  [font=\tiny,xslant=0.02] [align=left] {$\displaystyle \{2\} \leq \emptyset $};
\draw (401.87,112.33) node [anchor=north west][inner sep=0.75pt]  [font=\tiny,xslant=0.02] [align=left] {$\displaystyle \{1\} \leq \emptyset $};
\draw (344.56,171.47) node [anchor=north west][inner sep=0.75pt]  [font=\tiny,xslant=0.02] [align=left] {$\displaystyle \emptyset \leq \{\{1\}\}$};
\draw (256.27,113.33) node [anchor=north west][inner sep=0.75pt]  [font=\tiny,xslant=0.02] [align=left] {$\displaystyle \emptyset \leq \{\{2\}\}$};
\draw (260.95,166.96) node [anchor=north west][inner sep=0.75pt]  [font=\tiny,xslant=0.02] [align=left] {$\displaystyle \emptyset \leq \{\emptyset \}$};
\draw (247.97,139.71) node [anchor=north west][inner sep=0.75pt]  [font=\tiny,xslant=0.02] [align=left] {$\displaystyle \emptyset \leq \{\emptyset ,\{2\}\}$};
\draw (300.81,180.97) node [anchor=north west][inner sep=0.75pt]  [font=\tiny,xslant=0.02] [align=left] {$\displaystyle \emptyset \leq \{\emptyset ,\{1\}\}$};
\draw (365.56,79.51) node [anchor=north west][inner sep=0.75pt]  [font=\tiny,xslant=0.02] [align=left] {$\displaystyle \{1,2\} \leq \emptyset $};
\draw (270,85.36) node [anchor=north west][inner sep=0.75pt]  [font=\tiny,xslant=0.02] [align=left] {$\displaystyle \{2\} \leq \{\{2\}\}$};
\draw (374.41,146.11) node [anchor=north west][inner sep=0.75pt]  [font=\tiny,xslant=0.02] [align=left] {$\displaystyle \{1\} \leq \{\{1\}\}$};
\end{tikzpicture}
\end{center}
One can show $\widetilde{\Sigma}_{B_2}$ is the normal fan of $\widetilde{\Pi}_2$. Thus the augmented Bergman complex is $\partial\widetilde{\Pi}_2^*$.
\end{exa}
Note that the numerator in (\ref{AugChow}) is the Stanley-Reisner ring of $\widetilde{\Sigma}_M$, since for each $i\in [n]$, $F\in\L(M)-\{[n]\}$ the rays $\sigma_{\{i\}\le\emptyset}$ corresponds to $y_i$, and $\sigma_{\emptyset\le\{F\}}$ corresponds to $x_F$.

Now consider the case of $B_n$, the augmented Bergman fan $\widetilde{\Sigma}_{B_n}$ is $\widetilde{\Sigma}_n$. Some easy calculations show that all conditions in Lemma \ref{CohomologyFaceRing} hold. Therefore we have 
\[
	\widetilde{A}(B_n)=\frac{\mathbb{Q}[\partial\widetilde{\Pi}_n^*]}{\langle\theta_1,\ldots, \theta_n \rangle}\cong_{\S_n} H^*(X_{\widetilde{\Sigma}_n},\mathbb{Q}).
\]
However, there was no analogous Feichtner--Yuzvinsky's result known for the augmented Chow ring. In order to overcome this, we turn our focus on another way of constructing the stellohedron $\widetilde{\Pi}_n$~\textemdash~ as the graph associahedron of $n$-star graph $K_{1,n}$. 

\subsection{Stellohedron as the graph associahedron of $K_{1,n}$}\label{Sec:GraphAss}
Let $G=(V,E)$ be a simple graph, the \emph{graphical building set} $\B(G)$ is the set of nonempty subsets $I$ of $V$ such that the induced subgraph of $G$ on $I$ is connected. In fact, $\B(G)$ is a building set of the Boolean lattice over $V$.\\
The \emph{$n$-star graph} $K_{1,n}$ has vertex set $V=[n]\cup\{*\}$ and edge set $E=\{\{i,*\}:i\in [n]\}$.
\begin{exa}
$\mathcal{B}(K_{1,2})$ consists of the following elements:
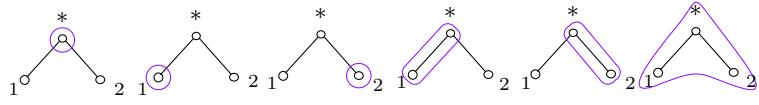
\begin{center}
\begin{tikzpicture}[x=0.6pt,y=0.6pt,yscale=-1,xscale=1]

\draw   (44.37,31.05) .. controls (44.28,29.78) and (45.34,28.69) .. (46.74,28.61) .. controls (48.14,28.52) and (49.35,29.49) .. (49.44,30.76) .. controls (49.53,32.02) and (48.47,33.12) .. (47.07,33.2) .. controls (45.67,33.29) and (44.46,32.32) .. (44.37,31.05) -- cycle ;
\draw   (20.27,57.26) .. controls (20.17,55.72) and (21.34,54.39) .. (22.88,54.29) .. controls (24.42,54.19) and (25.75,55.36) .. (25.85,56.9) .. controls (25.95,58.44) and (24.78,59.77) .. (23.24,59.87) .. controls (21.7,59.97) and (20.37,58.81) .. (20.27,57.26) -- cycle ;
\draw   (67.93,57.08) .. controls (67.83,55.73) and (69,54.57) .. (70.54,54.48) .. controls (72.08,54.4) and (73.41,55.42) .. (73.51,56.76) .. controls (73.61,58.11) and (72.44,59.27) .. (70.9,59.36) .. controls (69.36,59.45) and (68.03,58.43) .. (67.93,57.08) -- cycle ;
\draw    (45.12,32.57) -- (24.34,55.39) ;
\draw    (48.68,32.27) -- (69.03,55.24) ;
\draw   (128.7,29.46) .. controls (128.61,28.19) and (129.67,27.09) .. (131.07,27.01) .. controls (132.47,26.93) and (133.68,27.89) .. (133.77,29.16) .. controls (133.86,30.43) and (132.8,31.52) .. (131.4,31.61) .. controls (130,31.69) and (128.79,30.73) .. (128.7,29.46) -- cycle ;
\draw   (104.6,55.67) .. controls (104.5,54.13) and (105.67,52.8) .. (107.21,52.7) .. controls (108.75,52.6) and (110.08,53.76) .. (110.18,55.3) .. controls (110.28,56.85) and (109.11,58.18) .. (107.57,58.28) .. controls (106.03,58.38) and (104.7,57.21) .. (104.6,55.67) -- cycle ;
\draw   (152.26,55.48) .. controls (152.16,54.14) and (153.33,52.97) .. (154.87,52.89) .. controls (156.41,52.8) and (157.74,53.82) .. (157.84,55.17) .. controls (157.94,56.51) and (156.77,57.68) .. (155.23,57.76) .. controls (153.69,57.85) and (152.36,56.83) .. (152.26,55.48) -- cycle ;
\draw    (129.46,30.97) -- (108.67,53.79) ;
\draw    (133.02,30.67) -- (153.36,53.64) ;
\draw   (207.37,28.46) .. controls (207.28,27.19) and (208.34,26.09) .. (209.74,26.01) .. controls (211.14,25.93) and (212.35,26.89) .. (212.44,28.16) .. controls (212.53,29.43) and (211.47,30.52) .. (210.07,30.61) .. controls (208.67,30.69) and (207.46,29.73) .. (207.37,28.46) -- cycle ;
\draw   (183.27,54.67) .. controls (183.17,53.13) and (184.34,51.8) .. (185.88,51.7) .. controls (187.42,51.6) and (188.75,52.76) .. (188.85,54.3) .. controls (188.95,55.85) and (187.78,57.18) .. (186.24,57.28) .. controls (184.7,57.38) and (183.37,56.21) .. (183.27,54.67) -- cycle ;
\draw   (230.93,54.48) .. controls (230.83,53.14) and (232,51.97) .. (233.54,51.89) .. controls (235.08,51.8) and (236.41,52.82) .. (236.51,54.17) .. controls (236.61,55.51) and (235.44,56.68) .. (233.9,56.76) .. controls (232.36,56.85) and (231.03,55.83) .. (230.93,54.48) -- cycle ;
\draw    (208.12,29.97) -- (187.34,52.79) ;
\draw    (211.68,29.67) -- (232.03,52.64) ;
\draw   (289.03,27.79) .. controls (288.94,26.52) and (290,25.42) .. (291.41,25.34) .. controls (292.81,25.26) and (294.01,26.22) .. (294.11,27.49) .. controls (294.2,28.76) and (293.13,29.86) .. (291.73,29.94) .. controls (290.33,30.02) and (289.13,29.06) .. (289.03,27.79) -- cycle ;
\draw   (264.94,54) .. controls (264.84,52.46) and (266,51.13) .. (267.54,51.03) .. controls (269.08,50.93) and (270.41,52.1) .. (270.51,53.64) .. controls (270.61,55.18) and (269.45,56.51) .. (267.91,56.61) .. controls (266.37,56.71) and (265.04,55.54) .. (264.94,54) -- cycle ;
\draw   (312.6,53.82) .. controls (312.5,52.47) and (313.66,51.31) .. (315.21,51.22) .. controls (316.75,51.13) and (318.08,52.15) .. (318.18,53.5) .. controls (318.28,54.85) and (317.11,56.01) .. (315.57,56.1) .. controls (314.03,56.18) and (312.7,55.16) .. (312.6,53.82) -- cycle ;
\draw    (289.79,29.3) -- (269,52.13) ;
\draw    (293.35,29) -- (313.69,51.98) ;
\draw   (366.37,27.46) .. controls (366.28,26.19) and (367.34,25.09) .. (368.74,25.01) .. controls (370.14,24.93) and (371.35,25.89) .. (371.44,27.16) .. controls (371.53,28.43) and (370.47,29.52) .. (369.07,29.61) .. controls (367.67,29.69) and (366.46,28.73) .. (366.37,27.46) -- cycle ;
\draw   (342.27,53.67) .. controls (342.17,52.13) and (343.34,50.8) .. (344.88,50.7) .. controls (346.42,50.6) and (347.75,51.76) .. (347.85,53.3) .. controls (347.95,54.85) and (346.78,56.18) .. (345.24,56.28) .. controls (343.7,56.38) and (342.37,55.21) .. (342.27,53.67) -- cycle ;
\draw   (389.93,53.48) .. controls (389.83,52.14) and (391,50.97) .. (392.54,50.89) .. controls (394.08,50.8) and (395.41,51.82) .. (395.51,53.17) .. controls (395.61,54.51) and (394.44,55.68) .. (392.9,55.76) .. controls (391.36,55.85) and (390.03,54.83) .. (389.93,53.48) -- cycle ;
\draw    (367.12,28.97) -- (346.34,51.79) ;
\draw    (370.68,28.67) -- (391.03,51.64) ;
\draw  [color={rgb, 255:red, 144; green, 19; blue, 254 }  ,draw opacity=1 ] (39.14,31.84) .. controls (39.14,27.58) and (42.58,24.14) .. (46.83,24.14) .. controls (51.09,24.14) and (54.53,27.58) .. (54.53,31.84) .. controls (54.53,36.09) and (51.09,39.53) .. (46.83,39.53) .. controls (42.58,39.53) and (39.14,36.09) .. (39.14,31.84) -- cycle ;
\draw  [color={rgb, 255:red, 144; green, 19; blue, 254 }  ,draw opacity=1 ] (99.69,55.49) .. controls (99.69,51.23) and (103.14,47.79) .. (107.39,47.79) .. controls (111.64,47.79) and (115.09,51.23) .. (115.09,55.49) .. controls (115.09,59.74) and (111.64,63.18) .. (107.39,63.18) .. controls (103.14,63.18) and (99.69,59.74) .. (99.69,55.49) -- cycle ;
\draw  [color={rgb, 255:red, 144; green, 19; blue, 254 }  ,draw opacity=1 ] (226.02,54.33) .. controls (226.02,50.07) and (229.47,46.63) .. (233.72,46.63) .. controls (237.97,46.63) and (241.42,50.07) .. (241.42,54.33) .. controls (241.42,58.58) and (237.97,62.02) .. (233.72,62.02) .. controls (229.47,62.02) and (226.02,58.58) .. (226.02,54.33) -- cycle ;
\draw  [color={rgb, 255:red, 144; green, 19; blue, 254 }  ,draw opacity=1 ] (262.78,55.96) .. controls (260.73,54.07) and (260.59,50.88) .. (262.48,48.82) -- (286.07,23.18) .. controls (287.96,21.13) and (291.15,21) .. (293.21,22.89) -- (296.01,25.46) .. controls (298.07,27.35) and (298.2,30.55) .. (296.31,32.61) -- (272.73,58.24) .. controls (270.84,60.3) and (267.64,60.43) .. (265.58,58.54) -- cycle ;
\draw   (443.7,26.46) .. controls (443.61,25.19) and (444.67,24.09) .. (446.07,24.01) .. controls (447.47,23.93) and (448.68,24.89) .. (448.77,26.16) .. controls (448.86,27.43) and (447.8,28.52) .. (446.4,28.61) .. controls (445,28.69) and (443.79,27.73) .. (443.7,26.46) -- cycle ;
\draw   (419.6,52.67) .. controls (419.5,51.13) and (420.67,49.8) .. (422.21,49.7) .. controls (423.75,49.6) and (425.08,50.76) .. (425.18,52.3) .. controls (425.28,53.85) and (424.11,55.18) .. (422.57,55.28) .. controls (421.03,55.38) and (419.7,54.21) .. (419.6,52.67) -- cycle ;
\draw   (467.26,52.48) .. controls (467.16,51.14) and (468.33,49.97) .. (469.87,49.89) .. controls (471.41,49.8) and (472.74,50.82) .. (472.84,52.17) .. controls (472.94,53.51) and (471.77,54.68) .. (470.23,54.76) .. controls (468.69,54.85) and (467.36,53.83) .. (467.26,52.48) -- cycle ;
\draw    (444.46,27.97) -- (423.67,50.79) ;
\draw    (448.02,27.67) -- (468.36,50.64) ;
\draw  [color={rgb, 255:red, 144; green, 19; blue, 254 }  ,draw opacity=1 ] (367.89,21.7) .. controls (370.03,19.91) and (373.22,20.19) .. (375.01,22.33) -- (397.37,49.05) .. controls (399.16,51.19) and (398.88,54.38) .. (396.74,56.17) -- (393.81,58.61) .. controls (391.67,60.41) and (388.49,60.12) .. (386.69,57.98) -- (364.34,31.26) .. controls (362.55,29.12) and (362.83,25.94) .. (364.97,24.14) -- cycle ;
\draw  [color={rgb, 255:red, 144; green, 19; blue, 254 }  ,draw opacity=1 ] (465.93,31.4) .. controls (490.6,63.4) and (489.27,67.4) .. (468.6,60.07) .. controls (447.93,52.74) and (449.27,49.4) .. (422.6,60.07) .. controls (395.93,70.74) and (419.27,41.4) .. (434.6,22.74) .. controls (449.93,4.07) and (441.27,-0.6) .. (465.93,31.4) -- cycle ;

\draw (11.21,57.44) node [anchor=north west][inner sep=0.75pt]  [font=\tiny] [align=left] {$\displaystyle 1$};
\draw (41.27,12.6) node [anchor=north west][inner sep=0.75pt]  [font=\footnotesize] [align=left] {$\displaystyle \ast $};
\draw (77.99,57.26) node [anchor=north west][inner sep=0.75pt]  [font=\tiny] [align=left] {$\displaystyle 2$};
\draw (92.54,57.51) node [anchor=north west][inner sep=0.75pt]  [font=\tiny] [align=left] {$\displaystyle 1$};
\draw (125.6,11) node [anchor=north west][inner sep=0.75pt]  [font=\footnotesize] [align=left] {$\displaystyle \ast $};
\draw (162.33,52.66) node [anchor=north west][inner sep=0.75pt]  [font=\tiny] [align=left] {$\displaystyle 2$};
\draw (174.21,53.51) node [anchor=north west][inner sep=0.75pt]  [font=\tiny] [align=left] {$\displaystyle 1$};
\draw (204.27,10) node [anchor=north west][inner sep=0.75pt]  [font=\footnotesize] [align=left] {$\displaystyle \ast $};
\draw (240.99,54.66) node [anchor=north west][inner sep=0.75pt]  [font=\tiny] [align=left] {$\displaystyle 2$};
\draw (255.87,52.84) node [anchor=north west][inner sep=0.75pt]  [font=\tiny] [align=left] {$\displaystyle 1$};
\draw (285.93,9.33) node [anchor=north west][inner sep=0.75pt]  [font=\footnotesize] [align=left] {$\displaystyle \ast $};
\draw (320.66,54) node [anchor=north west][inner sep=0.75pt]  [font=\tiny] [align=left] {$\displaystyle 2$};
\draw (334.21,53.51) node [anchor=north west][inner sep=0.75pt]  [font=\tiny] [align=left] {$\displaystyle 1$};
\draw (363.27,9) node [anchor=north west][inner sep=0.75pt]  [font=\footnotesize] [align=left] {$\displaystyle \ast $};
\draw (399.99,53.66) node [anchor=north west][inner sep=0.75pt]  [font=\tiny] [align=left] {$\displaystyle 2$};
\draw (411.54,52.51) node [anchor=north west][inner sep=0.75pt]  [font=\tiny] [align=left] {$\displaystyle 1$};
\draw (440.6,8) node [anchor=north west][inner sep=0.75pt]  [font=\footnotesize] [align=left] {$\displaystyle \ast $};
\draw (476.33,52.66) node [anchor=north west][inner sep=0.75pt]  [font=\tiny] [align=left] {$\displaystyle 2$};
\end{tikzpicture}
\end{center}
\end{exa}
For graphical building sets $\B(G)$, the nested sets $N\subset\B(G)$ can be characterized by:
\begin{itemize}
\item[(1)] $\forall~I,J\in N$, either $I\subset J$, $I\supset J$ or $I\cap J=\emptyset$.
\item[(2)] $\forall~I,J\in N$ if $I\cap J=\emptyset$, then $I\cup J\notin\B(G)$ (not "connected").
\end{itemize}
From \cite[Theorem 6.5]{PRW2008}, $\widetilde{\N}(\L(B_{n+1}),\B(K_{1,n}))$ is combinatorially equivalent to $\partial\widetilde{\Pi}_n^*$. 
\begin{exa} \label{ExNestedComplexSt}
The reduced nested set complex with respect to $\B(K_{1,2})$ is  $\partial\widetilde{\Pi}_2^*$: 
\begin{center}
\begin{tikzpicture}[x=0.6pt,y=0.6pt,yscale=-1,xscale=1]

\draw   (81.68,202.38) .. controls (81.6,201.23) and (82.56,200.23) .. (83.83,200.16) .. controls (85.09,200.08) and (86.19,200.96) .. (86.27,202.11) .. controls (86.35,203.26) and (85.39,204.25) .. (84.12,204.33) .. controls (82.86,204.4) and (81.76,203.53) .. (81.68,202.38) -- cycle ;
\draw   (59.86,226.16) .. controls (59.77,224.76) and (60.83,223.55) .. (62.22,223.46) .. controls (63.62,223.37) and (64.82,224.43) .. (64.91,225.83) .. controls (65,227.22) and (63.94,228.43) .. (62.55,228.52) .. controls (61.15,228.61) and (59.95,227.55) .. (59.86,226.16) -- cycle ;
\draw   (103.01,225.99) .. controls (102.92,224.77) and (103.98,223.71) .. (105.38,223.63) .. controls (106.77,223.55) and (107.97,224.48) .. (108.07,225.7) .. controls (108.16,226.92) and (107.1,227.98) .. (105.7,228.06) .. controls (104.31,228.14) and (103.1,227.21) .. (103.01,225.99) -- cycle ;
\draw    (82.36,203.75) -- (63.54,224.45) ;
\draw    (85.59,203.48) -- (104.01,224.32) ;
\draw  [color={rgb, 255:red, 144; green, 19; blue, 254 }  ,draw opacity=1 ] (76.94,203.09) .. controls (76.94,199.23) and (80.06,196.1) .. (83.91,196.1) .. controls (87.76,196.1) and (90.88,199.23) .. (90.88,203.09) .. controls (90.88,206.95) and (87.76,210.07) .. (83.91,210.07) .. controls (80.06,210.07) and (76.94,206.95) .. (76.94,203.09) -- cycle ;

\draw   (134.2,18.89) .. controls (134.11,17.74) and (135.07,16.75) .. (136.34,16.67) .. controls (137.61,16.6) and (138.7,17.47) .. (138.79,18.62) .. controls (138.87,19.77) and (137.91,20.77) .. (136.64,20.84) .. controls (135.37,20.92) and (134.28,20.04) .. (134.2,18.89) -- cycle ;
\draw   (112.38,42.67) .. controls (112.29,41.27) and (113.34,40.06) .. (114.74,39.97) .. controls (116.13,39.88) and (117.34,40.94) .. (117.43,42.34) .. controls (117.52,43.74) and (116.46,44.94) .. (115.07,45.03) .. controls (113.67,45.13) and (112.47,44.07) .. (112.38,42.67) -- cycle ;
\draw   (155.53,42.5) .. controls (155.44,41.28) and (156.5,40.23) .. (157.89,40.15) .. controls (159.29,40.07) and (160.49,40.99) .. (160.58,42.21) .. controls (160.67,43.44) and (159.62,44.49) .. (158.22,44.57) .. controls (156.83,44.65) and (155.62,43.72) .. (155.53,42.5) -- cycle ;
\draw    (134.88,20.26) -- (116.06,40.97) ;
\draw    (138.1,19.99) -- (156.52,40.83) ;
\draw  [color={rgb, 255:red, 144; green, 19; blue, 254 }  ,draw opacity=1 ] (107.93,42.5) .. controls (107.93,38.65) and (111.05,35.52) .. (114.9,35.52) .. controls (118.75,35.52) and (121.87,38.65) .. (121.87,42.5) .. controls (121.87,46.36) and (118.75,49.49) .. (114.9,49.49) .. controls (111.05,49.49) and (107.93,46.36) .. (107.93,42.5) -- cycle ;

\draw   (233.49,88.74) .. controls (233.41,87.59) and (234.37,86.6) .. (235.64,86.52) .. controls (236.91,86.45) and (238,87.32) .. (238.09,88.47) .. controls (238.17,89.62) and (237.21,90.62) .. (235.94,90.69) .. controls (234.67,90.77) and (233.58,89.89) .. (233.49,88.74) -- cycle ;
\draw   (211.67,112.52) .. controls (211.58,111.12) and (212.64,109.92) .. (214.04,109.83) .. controls (215.43,109.74) and (216.63,110.79) .. (216.73,112.19) .. controls (216.82,113.59) and (215.76,114.8) .. (214.36,114.89) .. controls (212.97,114.98) and (211.76,113.92) .. (211.67,112.52) -- cycle ;
\draw   (254.83,112.35) .. controls (254.74,111.13) and (255.8,110.08) .. (257.19,110) .. controls (258.59,109.92) and (259.79,110.84) .. (259.88,112.07) .. controls (259.97,113.29) and (258.91,114.34) .. (257.52,114.42) .. controls (256.12,114.5) and (254.92,113.58) .. (254.83,112.35) -- cycle ;
\draw    (234.18,90.12) -- (215.36,110.82) ;
\draw    (237.4,89.84) -- (255.82,110.68) ;
\draw  [color={rgb, 255:red, 144; green, 19; blue, 254 }  ,draw opacity=1 ] (250.38,112.21) .. controls (250.38,108.35) and (253.5,105.23) .. (257.35,105.23) .. controls (261.2,105.23) and (264.32,108.35) .. (264.32,112.21) .. controls (264.32,116.07) and (261.2,119.19) .. (257.35,119.19) .. controls (253.5,119.19) and (250.38,116.07) .. (250.38,112.21) -- cycle ;

\draw   (37.61,88.74) .. controls (37.53,87.59) and (38.49,86.6) .. (39.76,86.52) .. controls (41.03,86.45) and (42.12,87.32) .. (42.21,88.47) .. controls (42.29,89.62) and (41.33,90.62) .. (40.06,90.69) .. controls (38.79,90.77) and (37.7,89.89) .. (37.61,88.74) -- cycle ;
\draw   (15.79,112.52) .. controls (15.7,111.12) and (16.76,109.92) .. (18.16,109.83) .. controls (19.55,109.74) and (20.75,110.79) .. (20.85,112.19) .. controls (20.94,113.59) and (19.88,114.8) .. (18.48,114.89) .. controls (17.09,114.98) and (15.88,113.92) .. (15.79,112.52) -- cycle ;
\draw   (58.95,112.35) .. controls (58.86,111.13) and (59.92,110.08) .. (61.31,110) .. controls (62.71,109.92) and (63.91,110.84) .. (64,112.07) .. controls (64.09,113.29) and (63.03,114.34) .. (61.64,114.42) .. controls (60.24,114.5) and (59.04,113.58) .. (58.95,112.35) -- cycle ;
\draw    (38.3,90.12) -- (19.47,110.82) ;
\draw    (41.52,89.84) -- (59.94,110.68) ;
\draw  [color={rgb, 255:red, 144; green, 19; blue, 254 }  ,draw opacity=1 ] (13.84,114.3) .. controls (11.98,112.59) and (11.86,109.69) .. (13.57,107.82) -- (34.93,84.57) .. controls (36.64,82.7) and (39.53,82.58) .. (41.39,84.3) -- (43.93,86.63) .. controls (45.79,88.35) and (45.91,91.25) .. (44.2,93.11) -- (22.85,116.37) .. controls (21.13,118.23) and (18.24,118.35) .. (16.38,116.64) -- cycle ;

\draw   (198.78,202.14) .. controls (198.7,200.99) and (199.66,199.99) .. (200.93,199.92) .. controls (202.2,199.84) and (203.29,200.72) .. (203.38,201.87) .. controls (203.46,203.02) and (202.5,204.01) .. (201.23,204.09) .. controls (199.96,204.16) and (198.87,203.29) .. (198.78,202.14) -- cycle ;
\draw   (176.96,225.92) .. controls (176.87,224.52) and (177.93,223.31) .. (179.33,223.22) .. controls (180.72,223.13) and (181.93,224.19) .. (182.02,225.59) .. controls (182.11,226.99) and (181.05,228.19) .. (179.65,228.28) .. controls (178.26,228.37) and (177.06,227.31) .. (176.96,225.92) -- cycle ;
\draw   (220.12,225.75) .. controls (220.03,224.53) and (221.09,223.47) .. (222.48,223.39) .. controls (223.88,223.31) and (225.08,224.24) .. (225.17,225.46) .. controls (225.26,226.68) and (224.2,227.74) .. (222.81,227.82) .. controls (221.41,227.9) and (220.21,226.97) .. (220.12,225.75) -- cycle ;
\draw    (199.47,203.51) -- (180.65,224.22) ;
\draw    (202.69,203.24) -- (221.11,224.08) ;
\draw  [color={rgb, 255:red, 144; green, 19; blue, 254 }  ,draw opacity=1 ] (200.17,196.92) .. controls (202.1,195.29) and (204.99,195.55) .. (206.61,197.49) -- (226.86,221.73) .. controls (228.48,223.67) and (228.22,226.56) .. (226.28,228.19) -- (223.64,230.4) .. controls (221.7,232.03) and (218.81,231.77) .. (217.19,229.83) -- (196.95,205.59) .. controls (195.32,203.65) and (195.58,200.76) .. (197.52,199.13) -- cycle ;

\draw    (55.82,86.89) -- (101.09,47.58) ;
\draw    (169.9,45.16) -- (225.74,91.12) ;
\draw    (52.2,119.55) -- (66.08,189.7) ;
\draw    (233.29,117.73) -- (212.76,196.35) ;
\draw    (105.32,213.89) -- (172.32,214.5) ;
\draw   (247.35,151.77) .. controls (247.29,150.92) and (248.03,150.18) .. (249,150.13) .. controls (249.97,150.07) and (250.81,150.72) .. (250.87,151.57) .. controls (250.93,152.43) and (250.2,153.17) .. (249.23,153.22) .. controls (248.26,153.28) and (247.42,152.63) .. (247.35,151.77) -- cycle ;
\draw   (230.65,169.42) .. controls (230.58,168.39) and (231.39,167.49) .. (232.45,167.42) .. controls (233.52,167.35) and (234.44,168.14) .. (234.51,169.18) .. controls (234.58,170.22) and (233.77,171.11) .. (232.71,171.18) .. controls (231.64,171.25) and (230.72,170.46) .. (230.65,169.42) -- cycle ;
\draw   (263.69,169.3) .. controls (263.62,168.39) and (264.43,167.61) .. (265.5,167.55) .. controls (266.57,167.49) and (267.49,168.18) .. (267.56,169.09) .. controls (267.63,169.99) and (266.82,170.78) .. (265.75,170.83) .. controls (264.68,170.89) and (263.76,170.21) .. (263.69,169.3) -- cycle ;
\draw    (247.88,152.79) -- (233.46,168.16) ;
\draw    (250.35,152.59) -- (264.45,168.06) ;
\draw  [color={rgb, 255:red, 144; green, 19; blue, 254 }  ,draw opacity=1 ] (248.22,147.49) .. controls (249.81,146.2) and (252.18,146.41) .. (253.52,147.96) -- (271.54,168.87) .. controls (272.88,170.42) and (272.67,172.73) .. (271.08,174.02) -- (268.91,175.79) .. controls (267.32,177.08) and (264.94,176.87) .. (263.61,175.32) -- (245.59,154.4) .. controls (244.25,152.85) and (244.46,150.54) .. (246.05,149.25) -- cycle ;
\draw  [color={rgb, 255:red, 144; green, 19; blue, 254 }  ,draw opacity=1 ] (261.18,169.19) .. controls (261.18,166.73) and (263.17,164.74) .. (265.63,164.74) .. controls (268.08,164.74) and (270.07,166.73) .. (270.07,169.19) .. controls (270.07,171.65) and (268.08,173.65) .. (265.63,173.65) .. controls (263.17,173.65) and (261.18,171.65) .. (261.18,169.19) -- cycle ;

\draw   (215.47,38.79) .. controls (215.41,37.96) and (216.1,37.25) .. (217.02,37.19) .. controls (217.93,37.14) and (218.72,37.77) .. (218.78,38.6) .. controls (218.84,39.43) and (218.15,40.15) .. (217.23,40.2) .. controls (216.32,40.26) and (215.53,39.63) .. (215.47,38.79) -- cycle ;
\draw   (199.72,55.95) .. controls (199.65,54.95) and (200.42,54.07) .. (201.42,54.01) .. controls (202.43,53.94) and (203.3,54.71) .. (203.37,55.72) .. controls (203.43,56.73) and (202.67,57.6) .. (201.66,57.66) .. controls (200.65,57.73) and (199.79,56.96) .. (199.72,55.95) -- cycle ;
\draw   (230.86,55.83) .. controls (230.8,54.95) and (231.56,54.19) .. (232.57,54.13) .. controls (233.58,54.08) and (234.44,54.74) .. (234.51,55.63) .. controls (234.58,56.51) and (233.81,57.27) .. (232.81,57.33) .. controls (231.8,57.38) and (230.93,56.72) .. (230.86,55.83) -- cycle ;
\draw    (215.96,39.78) -- (202.38,54.73) ;
\draw    (218.29,39.59) -- (231.58,54.63) ;
\draw  [color={rgb, 255:red, 144; green, 19; blue, 254 }  ,draw opacity=1 ] (227.66,55.73) .. controls (227.66,52.95) and (229.91,50.69) .. (232.69,50.69) .. controls (235.47,50.69) and (237.72,52.95) .. (237.72,55.73) .. controls (237.72,58.51) and (235.47,60.77) .. (232.69,60.77) .. controls (229.91,60.77) and (227.66,58.51) .. (227.66,55.73) -- cycle ;
\draw  [color={rgb, 255:red, 144; green, 19; blue, 254 }  ,draw opacity=1 ] (196.51,55.84) .. controls (196.51,53.05) and (198.76,50.8) .. (201.54,50.8) .. controls (204.32,50.8) and (206.57,53.05) .. (206.57,55.84) .. controls (206.57,58.62) and (204.32,60.88) .. (201.54,60.88) .. controls (198.76,60.88) and (196.51,58.62) .. (196.51,55.84) -- cycle ;

\draw   (136.85,227.38) .. controls (136.78,226.56) and (137.5,225.84) .. (138.44,225.79) .. controls (139.38,225.73) and (140.19,226.36) .. (140.25,227.19) .. controls (140.31,228.02) and (139.6,228.73) .. (138.66,228.79) .. controls (137.72,228.84) and (136.91,228.21) .. (136.85,227.38) -- cycle ;
\draw   (120.66,244.48) .. controls (120.6,243.47) and (121.38,242.6) .. (122.41,242.54) .. controls (123.45,242.47) and (124.34,243.24) .. (124.41,244.24) .. controls (124.48,245.24) and (123.69,246.11) .. (122.66,246.18) .. controls (121.62,246.24) and (120.73,245.48) .. (120.66,244.48) -- cycle ;
\draw   (152.67,244.36) .. controls (152.6,243.48) and (153.39,242.72) .. (154.42,242.66) .. controls (155.45,242.61) and (156.35,243.27) .. (156.41,244.15) .. controls (156.48,245.03) and (155.7,245.79) .. (154.66,245.84) .. controls (153.63,245.9) and (152.73,245.23) .. (152.67,244.36) -- cycle ;
\draw    (137.35,228.37) -- (123.39,243.25) ;
\draw    (139.74,228.18) -- (153.4,243.16) ;
\draw  [color={rgb, 255:red, 144; green, 19; blue, 254 }  ,draw opacity=1 ] (135.52,222.18) .. controls (137,220.86) and (139.31,220.96) .. (140.68,222.4) -- (159.05,241.81) .. controls (160.42,243.25) and (160.32,245.49) .. (158.84,246.82) -- (156.82,248.62) .. controls (155.34,249.94) and (153.03,249.85) .. (151.67,248.41) -- (133.29,228.99) .. controls (131.93,227.55) and (132.02,225.31) .. (133.5,223.99) -- cycle ;
\draw  [color={rgb, 255:red, 144; green, 19; blue, 254 }  ,draw opacity=1 ] (134.42,227.38) .. controls (134.42,225) and (136.34,223.07) .. (138.72,223.07) .. controls (141.1,223.07) and (143.03,225) .. (143.03,227.38) .. controls (143.03,229.76) and (141.1,231.7) .. (138.72,231.7) .. controls (136.34,231.7) and (134.42,229.76) .. (134.42,227.38) -- cycle ;

\draw   (63.53,31.54) .. controls (63.47,30.65) and (64.21,29.88) .. (65.19,29.82) .. controls (66.17,29.76) and (67.02,30.44) .. (67.08,31.33) .. controls (67.15,32.22) and (66.4,32.99) .. (65.42,33.05) .. controls (64.44,33.1) and (63.59,32.43) .. (63.53,31.54) -- cycle ;
\draw   (46.65,49.94) .. controls (46.58,48.85) and (47.39,47.92) .. (48.47,47.85) .. controls (49.55,47.78) and (50.49,48.6) .. (50.56,49.68) .. controls (50.63,50.76) and (49.81,51.7) .. (48.73,51.77) .. controls (47.65,51.84) and (46.72,51.02) .. (46.65,49.94) -- cycle ;
\draw   (80.04,49.81) .. controls (79.97,48.86) and (80.79,48.04) .. (81.87,47.98) .. controls (82.95,47.92) and (83.88,48.64) .. (83.95,49.58) .. controls (84.02,50.53) and (83.2,51.35) .. (82.12,51.41) .. controls (81.04,51.47) and (80.11,50.75) .. (80.04,49.81) -- cycle ;
\draw    (64.06,32.6) -- (49.5,48.62) ;
\draw    (66.55,32.39) -- (80.81,48.51) ;
\draw  [color={rgb, 255:red, 144; green, 19; blue, 254 }  ,draw opacity=1 ] (43.55,53.45) .. controls (42.01,52.04) and (41.91,49.64) .. (43.33,48.11) -- (61.99,27.78) .. controls (63.41,26.24) and (65.79,26.14) .. (67.33,27.55) -- (69.42,29.48) .. controls (70.96,30.9) and (71.06,33.29) .. (69.65,34.83) -- (50.98,55.16) .. controls (49.57,56.7) and (47.18,56.8) .. (45.64,55.38) -- cycle ;
\draw  [color={rgb, 255:red, 144; green, 19; blue, 254 }  ,draw opacity=1 ] (44.55,49.95) .. controls (44.55,47.79) and (46.3,46.04) .. (48.46,46.04) .. controls (50.61,46.04) and (52.37,47.79) .. (52.37,49.95) .. controls (52.37,52.12) and (50.61,53.87) .. (48.46,53.87) .. controls (46.3,53.87) and (44.55,52.12) .. (44.55,49.95) -- cycle ;

\draw   (30.26,153.21) .. controls (30.19,152.34) and (30.92,151.59) .. (31.87,151.54) .. controls (32.82,151.48) and (33.64,152.14) .. (33.7,153) .. controls (33.77,153.87) and (33.04,154.61) .. (32.09,154.67) .. controls (31.14,154.73) and (30.32,154.07) .. (30.26,153.21) -- cycle ;
\draw   (13.87,171.06) .. controls (13.8,170.01) and (14.59,169.11) .. (15.64,169.04) .. controls (16.69,168.97) and (17.59,169.77) .. (17.66,170.82) .. controls (17.73,171.87) and (16.94,172.77) .. (15.89,172.84) .. controls (14.84,172.91) and (13.94,172.11) .. (13.87,171.06) -- cycle ;
\draw   (46.28,170.94) .. controls (46.21,170.02) and (47.01,169.23) .. (48.05,169.17) .. controls (49.1,169.11) and (50.01,169.8) .. (50.07,170.72) .. controls (50.14,171.64) and (49.35,172.43) .. (48.3,172.49) .. controls (47.25,172.55) and (46.35,171.86) .. (46.28,170.94) -- cycle ;
\draw    (30.77,154.24) -- (16.63,169.79) ;
\draw    (33.19,154.03) -- (47.02,169.68) ;
\draw  [color={rgb, 255:red, 144; green, 19; blue, 254 }  ,draw opacity=1 ] (11.54,173.8) .. controls (10.05,172.42) and (9.95,170.1) .. (11.32,168.61) -- (29.44,148.87) .. controls (30.81,147.38) and (33.13,147.28) .. (34.62,148.66) -- (36.66,150.53) .. controls (38.15,151.9) and (38.24,154.23) .. (36.87,155.72) -- (18.75,175.45) .. controls (17.38,176.95) and (15.06,177.04) .. (13.57,175.67) -- cycle ;
\draw  [color={rgb, 255:red, 144; green, 19; blue, 254 }  ,draw opacity=1 ] (28.15,153.37) .. controls (28.15,151.27) and (29.85,149.56) .. (31.94,149.56) .. controls (34.04,149.56) and (35.74,151.27) .. (35.74,153.37) .. controls (35.74,155.47) and (34.04,157.17) .. (31.94,157.17) .. controls (29.85,157.17) and (28.15,155.47) .. (28.15,153.37) -- cycle ;

\draw (169.19,225.45) node [anchor=north west][inner sep=0.75pt]  [font=\tiny] [align=left] {$\displaystyle 1$};
\draw (195.41,184.79) node [anchor=north west][inner sep=0.75pt]  [font=\footnotesize] [align=left] {$\displaystyle \ast $};
\draw (228.76,225.59) node [anchor=north west][inner sep=0.75pt]  [font=\tiny] [align=left] {$\displaystyle 2$};
\draw (7.11,111.15) node [anchor=north west][inner sep=0.75pt]  [font=\tiny] [align=left] {$\displaystyle 1$};
\draw (34.24,71.4) node [anchor=north west][inner sep=0.75pt]  [font=\footnotesize] [align=left] {$\displaystyle \ast $};
\draw (65.78,112.19) node [anchor=north west][inner sep=0.75pt]  [font=\tiny] [align=left] {$\displaystyle 2$};
\draw (202.99,111.15) node [anchor=north west][inner sep=0.75pt]  [font=\tiny] [align=left] {$\displaystyle 1$};
\draw (230.12,71.4) node [anchor=north west][inner sep=0.75pt]  [font=\footnotesize] [align=left] {$\displaystyle \ast $};
\draw (263.47,112.19) node [anchor=north west][inner sep=0.75pt]  [font=\tiny] [align=left] {$\displaystyle 2$};
\draw (100.98,44.02) node [anchor=north west][inner sep=0.75pt]  [font=\tiny] [align=left] {$\displaystyle 1$};
\draw (130.82,1.54) node [anchor=north west][inner sep=0.75pt]  [font=\footnotesize] [align=left] {$\displaystyle \ast $};
\draw (164.17,39.62) node [anchor=north west][inner sep=0.75pt]  [font=\tiny] [align=left] {$\displaystyle 2$};
\draw (51.18,225.99) node [anchor=north west][inner sep=0.75pt]  [font=\tiny] [align=left] {$\displaystyle 1$};
\draw (78.3,185.03) node [anchor=north west][inner sep=0.75pt]  [font=\footnotesize] [align=left] {$\displaystyle \ast $};
\draw (109.84,225.83) node [anchor=north west][inner sep=0.75pt]  [font=\tiny] [align=left] {$\displaystyle 2$};
\draw (223.52,168.17) node [anchor=north west][inner sep=0.75pt]  [font=\tiny] [align=left] {$\displaystyle 1$};
\draw (243.36,136.66) node [anchor=north west][inner sep=0.75pt]  [font=\footnotesize] [align=left] {$\displaystyle \ast $};
\draw (270.83,168.28) node [anchor=north west][inner sep=0.75pt]  [font=\tiny] [align=left] {$\displaystyle 2$};
\draw (190.76,53.99) node [anchor=north west][inner sep=0.75pt]  [font=\tiny] [align=left] {$\displaystyle 1$};
\draw (211.36,24.47) node [anchor=north west][inner sep=0.75pt]  [font=\footnotesize] [align=left] {$\displaystyle \ast $};
\draw (235.71,54.74) node [anchor=north west][inner sep=0.75pt]  [font=\tiny] [align=left] {$\displaystyle 2$};
\draw (113.61,243.16) node [anchor=north west][inner sep=0.75pt]  [font=\tiny] [align=left] {$\displaystyle 1$};
\draw (132.79,212.03) node [anchor=north west][inner sep=0.75pt]  [font=\footnotesize] [align=left] {$\displaystyle \ast $};
\draw (159.43,243.26) node [anchor=north west][inner sep=0.75pt]  [font=\tiny] [align=left] {$\displaystyle 2$};
\draw (36.7,48.08) node [anchor=north west][inner sep=0.75pt]  [font=\tiny] [align=left] {$\displaystyle 1$};
\draw (59.56,16.64) node [anchor=north west][inner sep=0.75pt]  [font=\footnotesize] [align=left] {$\displaystyle \ast $};
\draw (84.19,48.89) node [anchor=north west][inner sep=0.75pt]  [font=\tiny] [align=left] {$\displaystyle 2$};
\draw (4.06,169.16) node [anchor=north west][inner sep=0.75pt]  [font=\tiny] [align=left] {$\displaystyle 1$};
\draw (26.23,138.56) node [anchor=north west][inner sep=0.75pt]  [font=\footnotesize] [align=left] {$\displaystyle \ast $};
\draw (50.16,169.95) node [anchor=north west][inner sep=0.75pt]  [font=\tiny] [align=left] {$\displaystyle 2$};

\end{tikzpicture}
\end{center}
\end{exa}

Comparing Example \ref{ExNestedComplexSt} and Example \ref{ExAugBergB2} leads us to a combinatorial proof of the following result\footnote{This is a weaker result than \cite[Proposition 2.6]{BHMPW2020} which states that the augmented Bergman fan of $B_n$ is the normal fan of the stellohedron.}.
\begin{prop}\label{FaceStructrueBn}
The augmented Bergman fan (complex) of $B_n$ is combinatorially equivalent to the dual stellohedron $\widetilde{\Pi}_n^*$. Consequently, there is a poset isomorphism between their face lattices.
\end{prop}
The isomorphism is as follows. Each nested set represents a face of $\partial\widetilde{\Pi}_n^*$, and it corresponds to a compatible pair $I\le\F$ which represents a cone in $\widetilde{\Sigma}_{B_n}$.
\begin{exa}
For $n=6$, the following is a nested set with respect to $\B(K_{1,6})$ and the corresponding compatible pair.
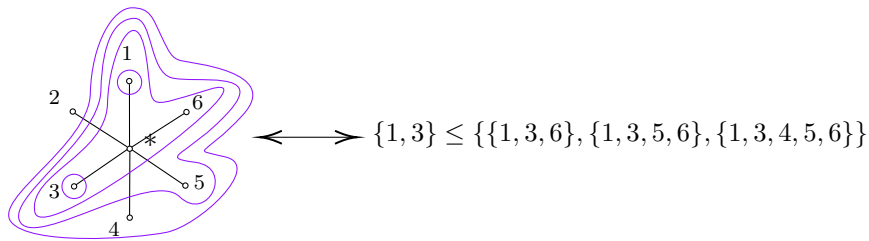
\[
\begin{tikzpicture}[x=0.75pt,y=0.75pt,yscale=-1,xscale=1]

\draw   (69.26,80.19) .. controls (69.22,79.45) and (69.77,78.81) .. (70.51,78.77) .. controls (71.25,78.72) and (71.88,79.28) .. (71.93,80.01) .. controls (71.98,80.75) and (71.42,81.38) .. (70.68,81.43) .. controls (69.95,81.48) and (69.31,80.92) .. (69.26,80.19) -- cycle ;
\draw    (71.86,80.73) -- (97.52,98.08) ;
\draw    (69.26,80.55) -- (43.52,98.3) ;
\draw    (70.51,78.77) -- (70.47,47.5) ;
\draw    (69.26,79.1) -- (43.06,62.29) ;
\draw    (71.57,79.29) -- (98.13,62.3) ;
\draw    (70.68,81.43) -- (70.53,113.56) ;
\draw   (97.89,61.66) .. controls (97.84,60.93) and (98.4,60.29) .. (99.13,60.24) .. controls (99.87,60.2) and (100.5,60.75) .. (100.55,61.49) .. controls (100.6,62.23) and (100.04,62.86) .. (99.31,62.91) .. controls (98.57,62.96) and (97.93,62.4) .. (97.89,61.66) -- cycle ;
\draw   (97.36,98.78) .. controls (97.31,98.05) and (97.87,97.41) .. (98.6,97.36) .. controls (99.34,97.31) and (99.98,97.87) .. (100.02,98.61) .. controls (100.07,99.34) and (99.51,99.98) .. (98.78,100.03) .. controls (98.04,100.08) and (97.41,99.52) .. (97.36,98.78) -- cycle ;
\draw   (40.55,61.41) .. controls (40.51,60.67) and (41.06,60.04) .. (41.8,59.99) .. controls (42.54,59.94) and (43.17,60.5) .. (43.22,61.23) .. controls (43.27,61.97) and (42.71,62.61) .. (41.97,62.65) .. controls (41.24,62.7) and (40.6,62.14) .. (40.55,61.41) -- cycle ;
\draw   (69.03,46.05) .. controls (68.98,45.32) and (69.54,44.68) .. (70.27,44.63) .. controls (71.01,44.59) and (71.64,45.14) .. (71.69,45.88) .. controls (71.74,46.62) and (71.18,47.25) .. (70.45,47.3) .. controls (69.71,47.35) and (69.08,46.79) .. (69.03,46.05) -- cycle ;
\draw   (69.28,114.98) .. controls (69.24,114.24) and (69.79,113.6) .. (70.53,113.56) .. controls (71.27,113.51) and (71.9,114.07) .. (71.95,114.8) .. controls (72,115.54) and (71.44,116.17) .. (70.7,116.22) .. controls (69.97,116.27) and (69.33,115.71) .. (69.28,114.98) -- cycle ;
\draw   (41.21,99.2) .. controls (41.17,98.46) and (41.72,97.83) .. (42.46,97.78) .. controls (43.2,97.73) and (43.83,98.29) .. (43.88,99.02) .. controls (43.93,99.76) and (43.37,100.4) .. (42.63,100.44) .. controls (41.9,100.49) and (41.26,99.93) .. (41.21,99.2) -- cycle ;
\draw  [color={rgb, 255:red, 144; green, 19; blue, 254 }  ,draw opacity=1 ] (64.4,46.54) .. controls (64.4,43.17) and (67.14,40.43) .. (70.52,40.43) .. controls (73.9,40.43) and (76.64,43.17) .. (76.64,46.54) .. controls (76.64,49.92) and (73.9,52.66) .. (70.52,52.66) .. controls (67.14,52.66) and (64.4,49.92) .. (64.4,46.54) -- cycle ;
\draw  [color={rgb, 255:red, 144; green, 19; blue, 254 }  ,draw opacity=1 ] (36.43,99.11) .. controls (36.43,95.73) and (39.17,92.99) .. (42.55,92.99) .. controls (45.93,92.99) and (48.66,95.73) .. (48.66,99.11) .. controls (48.66,102.49) and (45.93,105.23) .. (42.55,105.23) .. controls (39.17,105.23) and (36.43,102.49) .. (36.43,99.11) -- cycle ;
\draw  [color={rgb, 255:red, 144; green, 19; blue, 254 }  ,draw opacity=1 ] (59.4,57.3) .. controls (60.07,31.3) and (74.07,3.96) .. (78.73,32.63) .. controls (83.4,61.3) and (75.4,60.63) .. (111.4,49.96) .. controls (147.4,39.3) and (26.07,131.3) .. (22.73,111.96) .. controls (19.4,92.63) and (58.73,83.3) .. (59.4,57.3) -- cycle ;
\draw  [color={rgb, 255:red, 144; green, 19; blue, 254 }  ,draw opacity=1 ] (56.07,50.63) .. controls (56.73,24.63) and (63.73,14.3) .. (73.07,14.3) .. controls (82.4,14.3) and (80.3,41.26) .. (91.73,44.3) .. controls (103.17,47.33) and (140.4,39.63) .. (120.4,61.63) .. controls (100.4,83.63) and (83.73,84.3) .. (103.73,86.96) .. controls (123.73,89.63) and (109.89,118) .. (95.73,104.3) .. controls (81.58,90.6) and (25.07,129.63) .. (15.73,117.63) .. controls (6.4,105.63) and (55.4,76.63) .. (56.07,50.63) -- cycle ;
\draw  [color={rgb, 255:red, 144; green, 19; blue, 254 }  ,draw opacity=1 ] (50.73,52.3) .. controls (51.4,26.3) and (62.4,8.96) .. (71.73,8.96) .. controls (81.07,8.96) and (82.21,14.17) .. (86.4,22.96) .. controls (90.59,31.76) and (91.82,34.41) .. (96.4,35.63) .. controls (100.98,36.85) and (147.73,40.96) .. (127.73,62.96) .. controls (107.73,84.96) and (100.4,79.63) .. (119.07,84.3) .. controls (137.73,88.96) and (105.73,118.96) .. (89.07,122.3) .. controls (72.4,125.63) and (20.07,128.96) .. (10.73,116.96) .. controls (1.4,104.96) and (50.07,78.3) .. (50.73,52.3) -- cycle ;

\draw    (138.4,74.3) -- (183.07,74.3) ;
\draw [shift={(185.07,74.3)}, rotate = 180] [color={rgb, 255:red, 0; green, 0; blue, 0 }  ][line width=0.75]    (10.93,-3.29) .. controls (6.95,-1.4) and (3.31,-0.3) .. (0,0) .. controls (3.31,0.3) and (6.95,1.4) .. (10.93,3.29)   ;
\draw [shift={(136.4,74.3)}, rotate = 360] [color={rgb, 255:red, 0; green, 0; blue, 0 }  ][line width=0.75]    (10.93,-3.29) .. controls (6.95,-1.4) and (3.31,-0.3) .. (0,0) .. controls (3.31,0.3) and (6.95,1.4) .. (10.93,3.29)   ;

\draw (191.8,64.75) node [anchor=north west][inner sep=0.75pt]  [font=\small] [align=left] {$\displaystyle \{1,3\} \leq \{\{1,3,6\} ,\{1,3,5,6\} ,\{1,3,4,5,6\}\}$};
\draw (76.04,70.93) node [anchor=north west][inner sep=0.75pt]   [align=left] {$\displaystyle \ast $};
\draw (28.59,49.16) node [anchor=north west][inner sep=0.75pt]  [font=\scriptsize] [align=left] {$\displaystyle 2$};
\draw (100.83,51.76) node [anchor=north west][inner sep=0.75pt]  [font=\scriptsize] [align=left] {$\displaystyle 6$};
\draw (58.67,115.82) node [anchor=north west][inner sep=0.75pt]  [font=\scriptsize] [align=left] {$\displaystyle 4$};
\draw (65.12,27) node [anchor=north west][inner sep=0.75pt]  [font=\scriptsize] [align=left] {$\displaystyle 1$};
\draw (101.59,92.24) node [anchor=north west][inner sep=0.75pt]  [font=\scriptsize] [align=left] {$\displaystyle 5$};
\draw (28.59,97.06) node [anchor=north west][inner sep=0.75pt]  [font=\scriptsize] [align=left] {$\displaystyle 3$};

\end{tikzpicture}
\]
\end{exa}
We show in Section \ref{Sec:FaceAugBerg} that the isomorphism in Propostion \ref{FaceStructrueBn} holds even when the lattice is not a Boolean lattice. 

\subsection{Face Structure of augmented Bergman fan of $M$} \label{Sec:FaceAugBerg}
Let $M$, $\L(M)$, $\I(M)$ be defined as in Section \ref{Sec:AugChowBn}. The collection $\I(M)$ forms an abstract simplicial complex called the \emph{independence complex}; here we identify $I(M)$ with the face lattice of the independence complex. We construct a new poset $\widetilde{\mathcal{L}}(M)$ from $\mathcal{L}(M)$ and $\mathcal{I}(M)$ in the following ways:
\begin{itemize}
\item As a set, $\widetilde{\mathcal{L}}(M)=\mathcal{L}(M)\uplus\mathcal{I}(M)$. Denote by $F_{*}$ the flat $F$ in $\widetilde{\mathcal{L}}(M)$.
\item For $I\in\mathcal{I}(M)$, define $I\lessdot \mathsf{cl}_M(I)_*$ where $\mathsf{cl}_M(I)$ is the closure of $I$ in $M$.
The relations inside $\mathcal{L}(M)$ and $\mathcal{I}(M)$ stay the same.  
\end{itemize}

\begin{exa}
Consider the uniform matroid $U_{3,2}$, then the new poset is 
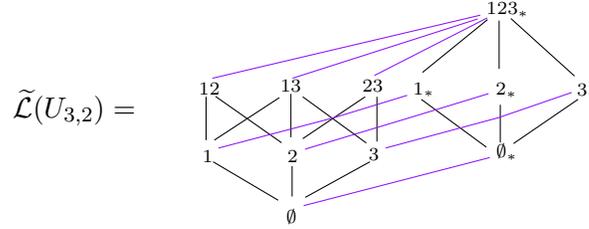
\begin{center}
\begin{tikzpicture}[x=0.75pt,y=0.75pt,yscale=-1.2,xscale=1.2]

\draw    (183.23,39.99) -- (210.56,17.77) ;
\draw    (215.42,18.19) -- (215.42,38.47) ;
\draw    (220.37,17.32) -- (247.11,41.5) ;
\draw    (182.73,51.32) -- (210.07,71.42) ;
\draw    (215.91,53.59) -- (215.91,69.45) ;
\draw    (221.46,69.91) -- (248.6,51.32) ;
\draw    (92.27,68.46) -- (92.27,50.23) ;
\draw    (164.08,69.07) -- (164.08,55.57) -- (164.08,50.23) ;
\draw    (96.28,50.11) -- (124.01,69.56) ;
\draw    (95.73,68.46) -- (123.07,50.59) ;
\draw    (127.92,49.01) -- (127.92,67.25) ;
\draw    (132.88,50.23) -- (159.62,69.68) ;
\draw    (132.88,67.86) -- (159.13,50.23) ;
\draw    (95.24,77.58) -- (122.58,93.75) ;
\draw    (128.42,79.4) -- (128.42,92.17) ;
\draw    (133.97,92.53) -- (161.11,77.58) ;
\draw [color={rgb, 255:red, 144; green, 19; blue, 254 }  ,draw opacity=1 ]   (211.51,75.79) -- (133.3,96.02) ;
\draw [color={rgb, 255:red, 144; green, 19; blue, 254 }  ,draw opacity=1 ]   (245.91,48.43) -- (215.12,59.03) -- (167.7,71.79) ;
\draw [color={rgb, 255:red, 144; green, 19; blue, 254 }  ,draw opacity=1 ]   (210.98,48.91) -- (180.19,58.32) -- (133.85,72.27) ;
\draw [color={rgb, 255:red, 144; green, 19; blue, 254 }  ,draw opacity=1 ]   (176.79,49.63) -- (146,59.04) -- (97.49,71.8) ;
\draw [color={rgb, 255:red, 144; green, 19; blue, 254 }  ,draw opacity=1 ]   (207.53,15.91) -- (94.92,42.44) ;
\draw [color={rgb, 255:red, 144; green, 19; blue, 254 }  ,draw opacity=1 ]   (206.44,17.5) -- (128.59,42.44) ;
\draw [color={rgb, 255:red, 144; green, 19; blue, 254 }  ,draw opacity=1 ]   (210.43,17.06) -- (162.99,41.29) ;

\draw (209.08,9.39) node [anchor=north west][inner sep=0.75pt]  [font=\tiny] [align=left] {$\displaystyle 123_{*}$};
\draw (178.26,43.13) node [anchor=north west][inner sep=0.75pt]  [font=\tiny] [align=left] {$\displaystyle 1_{*}$};
\draw (212.92,43.31) node [anchor=north west][inner sep=0.75pt]  [font=\tiny] [align=left] {$\displaystyle 2_{*}$};
\draw (247.08,43.61) node [anchor=north west][inner sep=0.75pt]  [font=\tiny] [align=left] {$\displaystyle 3_{*}$};
\draw (213.03,68.96) node [anchor=north west][inner sep=0.75pt]  [font=\tiny] [align=left] {$\displaystyle \emptyset_{*}$};
\draw (88.21,43.07) node [anchor=north west][inner sep=0.75pt]  [font=\tiny] [align=left] {$\displaystyle 12$};
\draw (122.41,42.13) node [anchor=north west][inner sep=0.75pt]  [font=\tiny] [align=left] {$\displaystyle 13$};
\draw (156.85,41.98) node [anchor=north west][inner sep=0.75pt]  [font=\tiny] [align=left] {$\displaystyle 23$};
\draw (89.68,71.75) node [anchor=north west][inner sep=0.75pt]  [font=\tiny] [align=left] {$\displaystyle 1$};
\draw (125.43,71.53) node [anchor=north west][inner sep=0.75pt]  [font=\tiny] [align=left] {$\displaystyle 2$};
\draw (159.59,70.66) node [anchor=north west][inner sep=0.75pt]  [font=\tiny] [align=left] {$\displaystyle 3$};
\draw (124.79,96.41) node [anchor=north west][inner sep=0.75pt]  [font=\tiny] [align=left] {$\displaystyle \emptyset $};
\draw (10,46) node [anchor=north west][inner sep=0.75pt]   [align=left] {$\widetilde{\mathcal{L}}(U_{3,2}) =$};

\end{tikzpicture}.
\end{center}
\end{exa}

\begin{lem}
Take $\mathcal{G}=\{\{1\},\ldots,\{n\}\}\cup\{F_{*}\}_{F\in\L(M)}$ as the building set in $\widetilde{\L}(M)$, then the nested sets are of the form 
\[
	\left\{\{i\}\right\}_{i\in I}\cup\left\{F_{*}\right\}_{F\in\mathcal{F}}
\]
for some compatible pair $I\le\mathcal{F}$ where $I\in\mathcal{I}(M)$ and $\mathcal{F}$ is a flag of arbitrary flats.
\end{lem}
With this lemma, we could recover $\widetilde{\Sigma}_M$ as the reduced nested set complex with respect to the building set $\G$.
\begin{thm}
$\widetilde{\N}(\widetilde{\mathcal{L}}(M),\mathcal{G})$ is combinatorially equivalent to the augmented Bergman fan (complex) of $M$. Consequently, there is a poset isomorphism between their face lattices:
\[
	\left\{\{i\}\right\}_{i\in I}\cup\left\{F_{*}\right\}_{F\in\mathcal{F}}\longleftrightarrow \sigma_{I\le\mathcal{F}}
\]
for compatible pair $I\le\mathcal{F}$ where $I\in\mathcal{I}(M)$ and flag $\mathcal{F}\subset\L(M)-\{[n]\}$.
\end{thm}
Furthermore, if we consider the Chow ring $D(\widetilde{\L}(M),\G)$, then
\begin{align*}
	D(\widetilde{\mathcal{L}}(M),\mathcal{G})
	&=\frac{\mathbb{Q}\left[\{x_F\}_{F\in\mathcal{F}(M)}\cup\{y_i\}_{i\in [n]}\right]/(I_1+I_2)}{\left\langle y_i+\sum_{F:i\in F}x_F\right\rangle_{i\in[n]}+\left\langle\sum_{F:\emptyset\subset F}x_F\right\rangle}
	=\frac{\mathbb{Q}\left[\{x_F\}_{F\in\mathcal{F}(M)\setminus [n]}\cup\{y_i\}_{i\in [n]}\right]/(I_1+I_2)}{\langle y_i-\sum_{F:i\notin F}x_F\rangle_{i\in[n]}}\\
	&=\widetilde{A}(M)
\end{align*}
Therefore, the augmented Chow ring of $M$ can be viewed as the Chow ring of $\widetilde{\L}(M)$ with respect to $\G$. We can apply Proposition \ref{FYbasis} to obtain a basis for $\widetilde{A}(M)$. 
\begin{cor}\label{AugBasis}
The augmented Chow ring $\widetilde{A}(M)$ of $M$ has the following basis
\[
	\widetilde{FY}(M)\coloneqq\left\{x_{F_1}^{a_1}x_{F_2}^{a_2}\ldots x_{F_k}^{a_k}: \substack{\emptyset=F_0\subsetneq F_1\subsetneq F_2\subsetneq\ldots\subsetneq F_k \\
    1\le a_1\le\rk(F_1),~ a_i\le\rk(F_i)-\rk(F_{i-1})-1 \text{ for }i\ge 2}
    \right\}
\]
\end{cor}

\begin{rem}
Mastroeni and Mccullough \cite{Mastroeni2022Koszul} mentioned that $\widetilde{A}(M)= D(\widetilde{\L}(M),\G)$ has also been discovered independently by Chris Eur, who further noticed that $\widetilde{\L}(M)$ is the lattice of flats of the \emph{free coextension} of $M$.
\end{rem}

\subsection{Back to Boolean matroids} \label{SecExtCodeBasis}
From Corollary \ref{AugBasis}, the basis $\widetilde{FY}(B_n)$ consists of monomials $x_{F_1}^{a_1}\ldots x_{F_k}^{a_k}$ indexed by a chain in the Boolean lattice with exponent $a_i$ satisfying
\begin{itemize}
	\item $|F_i|-|F_{i-1}|\ge 2$ for all $i\ge 2$. 
	\item $1\le a_1\le |F_1|$ and $a_i\le |F_i|-|F_{i-1}|-1$ for all $i\ge 2$.
\end{itemize}
\begin{exa}
$\widetilde{FY}(B_3)=\{1~\vline~x_1,x_2,x_3, x_{12}, x_{13}, x_{23},x_{123}~\vline~x_1x_{123},x_2x_{123},x_3x_{123},x_{12}^2,x_{13}^2,x_{23}^2,x_{123}^2~\vline~ x_{123}^3\}$.
\end{exa}
Similar to $FY(B_n)$, $\S_n$ also acts on $\widetilde{FY}(B_n)$ naturally and makes $\widetilde{A}(B_n)$ a permutation module. The following result is analogous to Theorem \ref{BasisBijection1}.
\begin{thm}\label{ExtendCodeBij}
There is a bijection $\widetilde{\phi}:\widetilde{FY}(B_n)\rightarrow \widetilde{\C}_n$ that respects the $\S_n$-actions and takes the degree of the monomials to the index$-1$ of the corresponding exteneded codes.
\end{thm}
The bijection is defined as follows. Given $u=x_{F_1}^{a_1}\ldots x_{F_k}^{a_k}\in \widetilde{FY}(B_n)$, let $\widetilde{\phi}(u)=(\alpha,f)$, where if $a_1=1$, then $\alpha_i=\begin{cases} j-1 & \text{ if }i\in F_j-F_{j-1}\\
\infty & \text{ if }i\in [n]-F_k
\end{cases}$
for $i\in [n]$, and $f(j)=a_{j+1}$ for $j\in [k]$; else if $a_1\ge 2$, then $\alpha_i=\begin{cases} j & \text{ if }i\in F_j-F_{j-1}\\
\infty & \text{ if }i\in [n]-F_k
\end{cases}$ for $i\in [n]$, and $f(1)=a_1-1$, $f(j)=a_j$ for $2\le j\le k$.
\begin{exa}
Let $u_1=x_{14}x_{1247}x_{1245679}^2$, $u_2=x_{14}^2x_{1247}x_{1245679}^2\in\widetilde{FY}(B_{9})$, then $\widetilde{\phi}(u_1)=01\infty022\hat{1}\infty\hat{2}$ and $\widetilde{\phi}(u_2)=12\infty\hat{1}33\hat{2}\infty\hat{3}$.
\end{exa}
Theorem \ref{ExtendCodeBij} shows that $\widetilde{A}^j(B_n)\otimes\mathbb{C}$ and $\widetilde{V}_{n,j-1}$ are isomorphic permutation representations of $\S_n$. This answers our question posed in the end of Section \ref{Sec:ExtCodes}. This also gives a new proof of Shareshian and Wachs' result (Theorem \ref{FrobStell}).

\section{Further work}
In this extended abstract, we mainly focus on the construction of FY-bases and their bijections with codes and extended codes. In the full paper  \cite{liao2022}, we also obtain enumerative results on codes, extended codes, and their connection with permutation statistics.

In another upcoming paper \cite{liao2022+}, we apply Corollary \ref{AugBasis} to uniform matroids $U_{n,r}$ and their $q$-analogs $M_r(\mathbb{F}_q^n)$, this enable us to obtain the counterpart of Hameister, Rao, and Simpson's  results on Chow ring in \cite{HameisterRaoSimpson2021} for the case of the augmented Chow rings. We further study the $\S_n$-module structure on $A(U_{n,r})$ and on $\widetilde{A}(U_{n,r})$, and obtain symmetric function analogs of the results mentioned above.

\section*{Acknowledgements}
The author thanks Michelle Wachs for her guidance and encouragement at every stage of this project and carefully reading the manuscript, as well as Vic Reiner for introducing the author to Chow rings of matroids, which motivated this work.

\bibliographystyle{abbrv}
\bibliography{bibliography}

\end{document}